\documentclass[reqno]{amsart}

\usepackage{lmodern}
\usepackage[T1]{fontenc}
\usepackage{textcomp}

\usepackage[utf8]{inputenc}
\usepackage[unicode]{hyperref}
\usepackage{amsfonts}
\usepackage{indentfirst}
\usepackage{enumerate}
\usepackage{amssymb}
\usepackage{amsmath}
\usepackage{mathrsfs}
\usepackage{paralist}
\usepackage{amsmath,amscd}
\usepackage{amsthm}

\DeclareMathOperator{\diam}{diam}
\DeclareMathOperator{\inte}{int}
\DeclareMathOperator{\dom}{dom}

\DeclareMathOperator{\cof}{cof}
\DeclareMathOperator{\cov}{cov}
\DeclareMathOperator{\non}{non}
\DeclareMathOperator{\Ker}{Ker}
\DeclareMathOperator{\length}{length}

\newtheorem{ttt}{Theorem}[section]
\newtheorem{llll}[ttt]{Lemma}
\newtheorem{ccc}[ttt]{Claim}
\newtheorem{eee}[ttt]{Example}
\newtheorem{fff}[ttt]{Fact}
\newtheorem{rrr}[ttt]{Remark}
\newtheorem{sss}[ttt]{Statement}
\newtheorem{ddd}[ttt]{Definition}
\newtheorem{qqq}[ttt]{Question}
\newtheorem{cccc}[ttt]{Corollary}
\newtheorem{nnn}[ttt]{Notation}
\newtheorem{ppp}[ttt]{Problem}
\newtheorem{pppp}[ttt]{Proposition}
\newtheorem{ccccc}[ttt]{Conjecture}

\newcommand{\beq}{\begin{equation} }

\newcommand{\bt}{\begin{ttt}}
\newcommand{\bl}{\begin{llll}}
\newcommand{\bc}{\begin{ccc}}
\newcommand{\bex}{\begin{eee}}
\newcommand{\bfa}{\begin{fff}}
\newcommand{\br}{\begin{rrr}\upshape}
\newcommand{\bst}{\begin{sss}}
\newcommand{\bd}{\begin{ddd}\upshape}
\newcommand{\bdd}{\begin{ddd}\upshape}

\newcommand{\bq}{\begin{qqq}}
\newcommand{\bnn}{\begin{nnn}}
\newcommand{\bpr}{\begin{ppp}}
\newcommand{\bprop}{\begin{pppp}}
\newcommand{\bcor}{\begin{cccc}}
\newcommand{\bcon}{\begin{ccccc}}

\newcommand{\eeq}{\end{equation}}
\newcommand{\et}{\end{ttt}}
\newcommand{\el}{\end{llll}}
\newcommand{\ec}{\end{ccc}}
\newcommand{\eex}{\end{eee}}
\newcommand{\efa}{\end{fff}}
\newcommand{\er}{\end{rrr}}
\newcommand{\est}{\end{sss}}
\newcommand{\ed}{\end{ddd}}
\newcommand{\eq}{\end{qqq}}
\newcommand{\ecor}{\end{cccc}}
\newcommand{\econ}{\end{ccccc}}
\newcommand{\enn}{\end{nnn}}
\newcommand{\epr}{\end{ppp}}
\newcommand{\eprop}{\end{pppp}}

\newcommand{\bp}{\noindent\textbf{Proof. }}
\newcommand{\ep}{\hspace{\stretch{1}}$\square$\medskip}

\begin{document}

\title{Answer to a question of Ros\l{}anowski and Shelah}
\author{Márk Poór}
\address{E\"otv\"os Lor\'and University, Institute of Mathematics, P\'azm\'any P\'eter
	s. 1/c, 1117 Budapest, Hungary}
\email{sokmark@cs.elte.hu}
\thanks{The author was supported by the National
	Research, Development and Innovation Office – NKFIH, grants no. 104178, 124749 and 129211}


\begin{abstract}
 In \cite{small-large} Ros\l{}anowski and Shelah asked whether every locally compact non-discrete group has a null but non-meager subgroup, 
 and conversely, whether it is consistent with $ZFC$ that in every locally compact group a meager subgroup is always null. They gave affirmative answers for both questions
 in the case of the Cantor group and the reals. In this paper we give affirmative answers for the general case.
\end{abstract}

\maketitle

\section{Introduction}
We are interested in null, but non-meager, and conversely meager but non-null subgroups of locally compact groups.
Only the non-discrete case is of interest, since nonempty open sets are of positive measure, thus in a discrete group $\{e \}$ has positive measure, and obviously is of second category.
It is known that under the Continuum Hypothesis one can construct null, but non-meager and meager, but non-null subgroups in the Cantor group and $\mathbb{R}$. (In fact the equalities $\cov(\mathcal{M}) = \cof(\mathcal{M})$ and  $\cov(\mathcal{N}) = \cof(\mathcal{N})$ are sufficient. Also, in general, it is easy to see that $\non(\mathcal{M})< \non(\mathcal{N})$ implies the existence of a null, non-meager subgroup in each locally compact Polish group, and similarly, from $\non(\mathcal{N})< \non(\mathcal{M})$ it follows that each locally compact Polish group admits a meager, but non-null subgroup.)
Talagrand proved in $ZFC$ that there exists a null, but non-meager filter in the Cantor group, that also yields a null, but non-meager subgroup \cite{talag}.
In \cite{small-large}, Ros\l{}anowski, and Shelah constructed $ZFC$ examples for null, but non-meager subgroups
of the Cantor group and the reals, and showed that it is consistent with $ZFC$ that in the two groups every meager
subgroup is null. Then they asked two questions

\textbf{Problem 5.1 \cite{small-large}}
\begin{enumerate}[(1)]
	\item  Does every non-discrete locally compact group (with complete Haar
	measure) admit a null non–meager subgroup?
	\item \label{mkerd} Is it consistent that no non-discrete locally compact group has a meager non–null subgroup?
\end{enumerate}
The fact, that the nonexistence of null, non-meager subgroups in the Cantor group is consistent was known earlier: 
H.M. Friedman proved that it is consistent with $ZFC$ that every $F_\sigma$ set in $2^\omega \times 2^\omega$ containing a rectangle of positive outer measure must contain a measurable rectangle of positive measure. Shelah and Fremlin proved that this statement implies the nonexistence of meager, non-null subgroups in $2^\omega$ (see \cite{Burke}).

In the first part of this paper we will construct a null but non-meager subgroup in
the case of the inverse limit of a countable sequence of finite groups, and in arbitrary second countable Lie-groups, and we show that 
these cases are sufficiently general.

In the second part, using Friedman's theorem we prove that it is consistent that in every locally compact Polish group every meager subgroup is null.
Last, by lemmas from the first part we reduce the case of arbitrary locally compact groups to the case of the Polish ones, thus conclude that the answer for $\eqref{mkerd}$ is affirmative too.

\section{Preliminaries and notations}


If we state that topological groups $G$ and $H$ are isomorphic, in symbols $G \simeq H$, then we mean that there is an algebraic isomorphism, which is also a homeomorphism.
Under the symbol $\leq$ we mean the subgroup relation, i.e. $H \leq G$ symbols that $H$ is a subgroup of $G$, and if we write
 $H \vartriangleleft G$, then $H$ is a normal subgroup of $G$. 
 All topological groups are assumed to be Hausdorff.
(It is known that for a locally compact Hausdorff topological group $G$, if $C \vartriangleleft G$ is a compact normal subgroup, then
$G/C$ is locally compact and Hausdorff.)

 A topological group  $G$ is locally compact if for each $g \in G$, there is a neighborhood $B$ of $g$, which is compact (i.e. each point has an open neighborhood which has compact closure).
Since in compact Hausdorff spaces the Baire category theorem is true (i.e. no nonempty open set can be covered by countably many nowhere dense sets), it is also true in locally compact Hausdorff spaces.

For any topological space $X$, $\mathcal{B}(X)$ denotes the Borel sets, i.e. the $\sigma$-algebra generated by the open sets.

Recall the following definition of left Haar Borel measure:
\bd \label{haardeff} Let $G$ be a locally compact group, and $\mu$ be a Borel measure, i.e.
	\[ \mu : \mathcal{B}(G) \to [0, \infty] \]
	such that
	\begin{itemize}
		\item	$\mu$ is left-invariant, i.e. for every $g \in G$, $B \in \mathcal{B}(G)$
			\[ \mu(gB) = \mu(B), \]
		\item $\mu(U) > 0$ for each open $U \neq \emptyset$,
		\item $\mu(K) < \infty$ for each compact set $K$,
		\item $\mu$ is inner regular with respect to the compact sets, that is for every  Borel $B$
		  \[  \mu(B) = \sup \{ \mu(K): \ K \subseteq B, \ K \text{ is compact} \} .\]
	\end{itemize}
	Then $\mu$ is called a left Haar Borel measure of $G$.
\ed
It is known that a left Haar Borel measure always exists, and it is unique up to a positive multiplicative constant.

Under left Haar measure we will mean the completion of a left Haar Borel measure.
\bd
	Let $G$ be a locally compact group. Suppose that $\nu$ is a left Haar Borel measure, and $\mu$ is the completion of $\nu$, that is
	a subset $H \subseteq G$ is measurable if it differs from a Borel set by (at most) a null Borel set, and then
	\[ \text{ if } (\exists B, B' \in \mathcal{B}(G):  B \bigtriangleup H \subseteq B' \ \wedge \nu(B') = 0) \ \textrm{ then } \mu(H):= \nu(B). \]
	Then $\mu$ is a left Haar measure of $G$.
\ed


From now on for a locally compact group $G$ $\mu_G$ will always symbol a left Haar measure.

It is easy to see (from the inner regularity) that if a measurable set $H$ is locally null, then it is null, i.e.
\begin{equation} \label{loknull}
 \left( \exists \text{ open cover }(U_\alpha)_{\alpha \in I} \text{ of } H  \ : \ (\forall \alpha)\  \mu_G(H \cap U_\alpha) = 0 \right) \ \Rightarrow \ \mu_G(H) = 0 .
\end{equation}

\br \upshape
	There is another definition of Haar Borel measure, namely, instead of inner regularity wrt. the compact sets we can require 
	\begin{itemize}
		\item outer regularity wrt. the open sets:
		\[  \mu(B) = \inf\{ \mu(U): \ B \subseteq U, \ U \text{ is open  }\} \text{ for every Borel } B,  \]
		\item inner regularity wrt. the compact sets, only in case of open sets:
		\[  \mu(U) = \sup \{ \mu(K): \ K \subseteq U, \ K \text{ is compact} \} \text{ for every } U \text{ open } \]
	\end{itemize}
	which also exists and is unique up to a positive multiplicative constant. 
	This definition results in a different left-invariant Borel measure only in the case of non-$\sigma$-compact groups.
	The proofs in the first chapter of this paper would also work with this definition of the Haar measure, but in the second 
	we will make use of the property $\eqref{loknull}$, and the theorem which we will present is simply not true if requiring these type of regularity properties, in Example $\ref{cex}$ we mention a well-known counterexample.
\er

\bd
	Let $(I, \preceq)$ be a directed partially ordered set (poset), i.e. for every $a,b \in I$ there exists a common upper bound $c$.
	And let $G_a$ ($a \in I$) be topological groups indexed by elements of $I$, and let $\varphi_{a,b}: G_b \to G_a$ for each pair $a \preceq b$ be surjective continuous homomorphisms, such that whenever $a \preceq b \preceq c$ then
	\[ \varphi_{a,c} = \varphi_{a,b} \circ \varphi_{b,c} \ . \]
	Then this system is called an inverse system, and the inverse limit of the system $\left( (G_a)_{a \in I}, (\varphi_{a,b})_{a \preceq b} \right)$ is
	\[ \underleftarrow{\lim}_{i \in I} G_i = \left\{  (x_i)_{i \in I} \in \prod_{i \in I} G_i : \ x_a = \varphi_{a,b}(x_b) \ (\forall a \preceq b \in I) \right\} \]
	with the topology inherited from the product space, and with the obvious (pointwise) group structure.
	
	If $I = \{ i \}$, then we identify $\underleftarrow{\lim}_{i \in I} G_i$ with $G_i$.

\ed

\br \upshape \label{bazis}
	By the definition of the product and subspace topology, there is a sub-base consisting of sets of the form
	\[ \left( \underleftarrow{\lim}_{a \in I} G_a \right) \cap \left( U_a \times \prod_{b \in I \setminus \{a \} } G_b \right) \ \ (a \in I, \ U_a \subseteq G_a \text{ open}) .\]
	But if $a \preceq d$ (since $(g_i)_{i \in I} \in \underleftarrow{\lim}_{a \in I} G_a \ \text{ implies } \ \varphi_{a,d}(g_d)=g_a$):
	\[ \left( \underleftarrow{\lim}_{i \in I} G_i \right) \cap \left( U_a \times \prod_{b \in I \setminus \{a \} } G_b \right) =
	\left( \underleftarrow{\lim}_{i \in I} G_i \right) \cap \left( \varphi_{a,d}^{-1}(U_a) \times \prod_{b \in I \setminus \{d \} } G_b \right). \]
	Thus these sets are closed under finite intersection (assume that $a,c \preceq d$), indeed,
	\[  \left( \left( \underleftarrow{\lim}_{i \in I} G_i \right) \cap \left( U_a \times \prod_{b \in I \setminus \{a \} } G_b \right) \right) \cap \left( \left(  \underleftarrow{\lim}_{i \in I} G_i \right) \cap \left( U_c \times \prod_{b \in I \setminus \{c \} } G_b \right) \right) =  \]
	\[ \left( \left( \underleftarrow{\lim}_{i \in I} G_i \right) \cap \left( \varphi_{a,d}^{-1}(U_a) \times \prod_{b \in I \setminus \{d \} } G_b \right) \right) \cap \left(  \left( \underleftarrow{\lim}_{i \in I} G_i \right) \cap \left( \varphi_{c,d}^{-1}(U_c) \times \prod_{b \in I \setminus \{d \} } G_b \right) \right)   = \]
	\[ = \left( \underleftarrow{\lim}_{i \in I} G_i \right) \cap \left( \left( \varphi_{c,d}^{-1}(U_c) \cap \varphi_{a,d}^{-1}(U_a) \right) \times \prod_{b \in I \setminus \{d \} } G_b \right). \]
	Hence, in fact, this sub-base is a base.
\er
It is known that for any locally compact group the null-ideals of a left Haar measure and a right Haar measure coincide, thus we can 
speak about null sets (\cite[442F ]{Fremlin}).
From now on $\mathcal{N}$ will denote the null-ideal, i.e.
\[ \mathcal{N} = \{ H \subseteq G:  \mu_G(H) = 0 \}. \]

\vspace{1cm}

\section{Null, but non-meager subgroups}

 Using ideas from \cite{small-large}, we will construct a null but non-meager subgroup in
inverse limits of countable many finite groups, and in  second countable Lie-groups.
First, from Lemma $\ref{gl}$ to Lemma $\ref{halmos}$, we will state lemmas so that we can restrict ourselves to those sufficiently general cases. We will find an appropriate  Lie group or profinite quotient for arbitrary locally compact groups, and we will show that the pull-back of a null but non-meager subgroup in that quotient is null and non-meager in the original group.
\br \upshape
 Using  ideas from \cite{elektoth}, in the case of locally compact Polish groups we could much more easily reduce the general case to the case of Lie groups and profinite groups. Moreover, if considering only locally compact Abelian Polish groups, then we only would have to construct null but non-meager subgroups in $\mathbb{R}$, countable product of finite groups and the $p$-adic integers.
\er

After that, by Proposition $\ref{profinite}$ the profinite case is handled, and later there will be a
similar construction for the Lie groups (Proposition $\ref{Lieeset}$). Finally we summarize, and state our main result, Theorem $\ref{atetel}$.

Recall the following lemma \cite[Thm. 1.6.1 ]{Tao}. 
\bl \label{gl}
	(Gleason-Yamabe theorem, stronger version). Let G be a
	locally compact group. Then there exists an open subgroup $G'$ of $G$ such that,
	for any open neighborhood $U$ of the identity in $G'$, there exists a compact
	normal subgroup $L_U \subseteq U$ of $G'$
	 such that $G'/L_U$ is a Lie group.
\el 

We will also make use of the following lemma \cite[Lemma 1.6.]{gl2}.
\bl \label{g1}
	Let $G$ be a locally compact group. Suppose that $N_1,N_2, \dots, N_k$ are closed normal subgroups, such that each for
	each $i$, $G/N_i$ is a Lie group. Then $G/(N_1 \cap N_2 \cap \dots \cap N_k)$ is a Lie group.
\el

From now on unless otherwise stated $\varphi_C$ always denotes the canonical projection from $G$ onto $G/C$, where $C$ is a compact normal subgroup.

\bd
	Suppose that $G$ is a topological group. Then $G_e$ denotes the connected component of the identity.
\ed

\bl \label{lie2}
	Let $G$ be a locally compact group, and $N,L \vartriangleleft G$ be compact normal subgroups, such that $G/N$, $G/L$ are Lie groups, and  assume that $G/N$ has $\lambda$-many connected components, i.e. $|G/N \ : \ (G/N)_e| = \lambda$. Then
	there is a finite number $k < \omega$, such that  $G/(L\cap N)$ has $\leq k \lambda$-many connected components.
\el

\bp
By Lemma $\ref{g1}$, $G/(L\cap N)$ is a Lie group.
	Now, applying the second isomorphism theorem for topological groups (see \cite[Thm. 1.5.18.]{arkh}),
	 \[ \left(G/(L \cap N)\right)/ \left(L/ (L \cap N) \right) \simeq G/L, \] and 
	\[ \left(G/(L \cap N)\right)/ \left(N/ (L \cap N) \right) \simeq G/N \] are Lie groups, thus replacing $G$ by $G/(L \cap N)$, (and $N$,$L$ by $N/(L \cap N)$, $L/(L \cap N)$, respectively)
   we can assume that $L \cap N= \{e \}$.
   
 Since	$N$ is a compact subgroup, by Cartan's theorem \cite[Thm. 20.12]{Lee} it is a   Lie group too. 
By its compactness, $N$ has only 
finitely many (say, $n$) connected components. We define $k=n$, and we will show that the open subgroup 
\begin{equation} \label{ekomp}
H:=\varphi_N^{-1}((G/N)_e)\leq G \textrm{  has at most } k \textrm{ connected components.}
\end{equation}
($N \vartriangleleft H$, and $\varphi_N|_H : H \to (G/N)_e$ is a surjective open continuous homomorphism, thus $H/N \simeq (G/N)_e$.)
Provided that $\eqref{ekomp}$ is true, it can be seen from the following that we are done. 
For any  surjective homomorphism $\psi: G_1 \to G_2$  between the groups $G_1$, $G_2$ and the  subgroup $X \leq G_2$
\[ |G_2 : X| = |G_1 : \psi^{-1}(X)|, \]
since the preimage of the disjoint cosets $gX$, $g'X$ are the disjoint sets $\psi^{-1}(gX) = \psi^{-1}(g) \psi^{-1}(X)$, 
$\psi^{-1}(g'X) = \psi^{-1}(g')\psi^{-1}(X)$, which are cosets of $\psi^{-1}(X)$. Thus
the index of $H$
\[ |G:H| = |\varphi_N^{-1}(G/N): \varphi_N^{-1}((G/N)_e)| = |G/N : (G/N)_e | = \lambda. \]
From which it clearly follows that the index $|G: H_e| = \lambda \cdot k$. But $H \leq G$ is an open subgroup of a Lie group, therefore $H$ is the union of some cosets of $G_e$, and $G_e = H_e$, thus $|G: G_e| = \lambda \cdot k$.

Turning to the proof of $\eqref{ekomp}$, let $\varphi_N'$ denote  the restriction of $\varphi_N$ to $H$, 
i.e. $\varphi_N' = \varphi_N|_H : H \to (G/N)_e$ is a surjective open continuous homomorphism.
Recall that $H \leq G$ is open, thus $H$ is a Lie group itself, thus $H_e$ is open in $H$.

It is obvious that $N_eH_e $ is a connected open subset of $H$ containing $H_e$
(continuous image of the connected set $N_e \times H_e $). 
Therefore it is the component of the identity, i.e. $N_eH_e = H_e$. 
This implies that if $N = \bigcup_{i<k} n_iN_e$, then 
\begin{equation} \label{szamolas}
 \bigcup_{i <k}n_iH_e = \bigcup_{i <k}n_iN_eH_e = \left( \bigcup_{i <k}n_iN_e \right) H_e = NH_e =  (\varphi_N')^{-1}(\varphi_N'(H_e)) 
\end{equation}

and that is a clopen subset of $H$ (the union of left cosets of the open subgroup $H_e$), so is its complement.
 But $NH_e$ is the preimage of the  set $\varphi_N(H_e) = \varphi_N(NH_e)$ under $\varphi_N$,
thus its complement
\[ H \setminus NH_e = \left( \varphi_N'^{-1}(H/N) \right) \setminus \left(  \varphi_N'^{-1}(\varphi'_N(NH_e)) \right)= \varphi_N'^{-1}\left( (H/N) \setminus \left(\varphi_N(NH_e)\right)\right) ) \]
is the preimage of $(H/N) \setminus (\varphi'_N(H_e))$. Therefore the clopen sets $NH_e$ and its complement is mapped by $\varphi_N'$
 onto the disjoint open sets $\varphi'_N(H_e) \subseteq H/N$ and its complement. But $H/N \simeq (G/N)_e$ is connected, thus at least one of these images is empty. Now $H_e \neq \emptyset$, thus $\varphi'_N(NH_e) \neq \emptyset$, which means that $H \setminus NH_e = \emptyset$, i.e. $H=NH_e$.
From this $H = NH_e = \bigcup_{i <k}n_iH_e$, i.e. $H$ is the union of (at most) $k$ connected open subsets.
\ep

\bd
	We call a locally compact group $G$ an FL-group, if for every  neighborhood $U \subseteq G$ of the identity, there is a  compact normal subgroup $N_U \subseteq U$ of $G$, s.t. $G/N_U$ is a Lie group that has finitely many connected components.
\ed

Now we can deduce the following.
\bl \label{struct}
	Let $G$ be a locally compact group.  Then there is an open FL subgroup $H$. 

\el
\bp
	Choose an arbitrary $G' \leq G$ which is given by Lemma $\ref{gl}$, and let $K \vartriangleleft G'$ be a compact normal subgroup, such that $G'/K$ is a Lie group.
	Let $\varphi_K: G' \to G'/K$ denote the canonical projection. Then, since $G'/K$ is a Lie group, $(G'/K)_e$ is open, thus
	$H:=\varphi_K^{-1}((G'/K)_e)$ is open in $G'$. But $G'$ is open in $G$, hence $H$ is open in $G$.
	
	First, we will check that $H$ satisfies the conclusion of Lemma $\ref{gl}$ (i.e. for every open set $U$ ($e \in U \subseteq H$) there exists a compact normal subgroup
	$L_U \vartriangleleft H$, such that $L_U \subseteq U$, and $H/L_U$ is a Lie group).
	
	Let $U \subseteq H$ be an open neighborhood of the identity. Then there is a compact normal subgroup $L_U \vartriangleleft G'$, such that $L_U \subseteq U$ and $G'/L_U$ is a Lie group.
	Then $L_U \vartriangleleft H \leq G'$ is obviously a normal subgroup of $H$. But if $\varphi_{L_U}: G' \to G'/L_U$ denotes the canonical projection, using that $H \subseteq G'$ is open, and $\varphi_{L_U}$ is an open mapping, $\varphi_{L_U}(H) = H/L_U \leq G'/L_U$ is an open subgroup of a Lie group, thus $H/L_U$ is a Lie group.
	
	Let $N_U = K \cap L_U$. Then $H/(K \cap L_U)$ is a Lie group by Lemma $\ref{g1}$, and since $H/K =  (G'/K)_e$ is connected, using Lemma $\ref{lie2}$ we obtain that
	$H/(K \cap L_U)$ has finitely many connected components. Thus choosing $N_U = K \cap L_U$ works.
\ep

Our plan is to take the subgroup $H$ given by the previous lemma, and for an appropriately choosen compact normal subgroup $L \vartriangleleft H$ we will construct a null but non-meager subgroup $S  \leq H/L$. 
Then at first we would like to ensure that the preimage of $S$ is also non-meager (in $H$), and hence it is also non-meager in $G$. That is what the following lemmas will be about.


We will need lemmas which are stated in terms of inverse limits.

For an arbitrary topological group $G$,
let $(\mathfrak{C}, \preceq)$ denote the directed system of compact normal subgroups, with the reversed inclusion relation, that is,
\[ \mathfrak{C} = \{N \vartriangleleft G : \ N \text{ is compact } \} ,\]
\[ M \preceq N \ \ \iff \ M \geq N ,\]

Now consider the following inverse system of topological groups.
\[  \left( \{ G/N : \ N \in \mathfrak{C} \}, \{ \varphi_{M,N}: G/N \to (G/N)/(M/N) \simeq G/M : \  M \preceq N \in \mathfrak{C}  \} \right). \]

\bd
	Let $G_a$ ($a \in I$) be an inverse system of topological groups ($(I, \preceq )$ is a directed set), and let $J \subseteq I$ be a directed subset (wrt. $\preceq$).
	Then let $\pi_{I \to J}: \underleftarrow{\lim}_{a \in I} G_a \to \underleftarrow{\lim}_{a \in J} G_a$ denote the natural projection
	\[ (x_a)_{a \in I} \mapsto (x_a)_{a \in J}, \]
	which is obviously a continuous homomorphism.
	Indeed, if we consider the projection on the whole product $\prod_{a \in I} G_a$, then it is continuous, and $\pi_{I \to J}$ is a restriction of it to the inverse limit.
	And clearly
	\[ \Ker(\pi_{I \to J}) = (\underleftarrow{\lim}_{a \in I} G_a) \cap \left( \prod_{a \in J} \{e_{G_a} \} \times \prod_{a \in I \setminus J} G_a \right). \].
	
\ed
	These projections may not be surjective.

Fix a locally compact FL-group $G$.
Let $(\mathfrak{N}, \preceq)$ denote the following subsystem of the system  $\mathfrak{C}$ (of compact normal subgroups of $G$)
\[ \mathfrak{N} = \{N \vartriangleleft G : \ N \in \mathfrak{C} \ \wedge \ G/N \text{ is a Lie group with finitely many connected components} \} \]
\[ N \preceq M \ \ \iff \ N \geq M.\]
Notice that, by Lemmas $\ref{g1}$,  and $\ref{lie2}$ we obtain the following lemma which we will make use of later. 

\bl \label{directed}
	Let $G$ be a locally compact group, and assume that there exists a compact normal subgroup $N \leq G$ such that $G/N$ is a Lie group with finitely many connected components. Then
 $\mathfrak{N}$ 
 is directed wrt. the reversed inclusion relation.
\el

\bl \label{reszrsz}
	Let $G$ be a topological group, and let $\mathfrak{M} \subseteq \mathfrak{C}$ be
	a directed subsystem of compact normal subgroups of $G$ (i.e. for each $M, M' \in \mathfrak{M}$ there is an $M'' \in \mathfrak{M}$ with
	$M,M'' \geq M''$), and set  $L:=\bigcap_{M \in \mathfrak{M}} M $. Then 
	\begin{itemize}
		\item the mapping
		\[ \psi: G \to  \underleftarrow{\lim}_{M \in \mathfrak{M}} G/M \]
		\[ g \mapsto 
		 (gM)_{M \in \mathfrak{M}} \]
		is a surjective open continuous homomorphism, and its kernel is $L$,
		\item moreover, the induced mapping
		$\overline{\psi}: G/L \to \underleftarrow{\lim}_{M \in \mathfrak{M}} G/M$ is an isomorphism.
	\end{itemize}
	
\el
\bp
First, $\psi$ is obviously a continuous homomorphism. 
 For the surjectivity, let $(g_\alpha M_\alpha)_{M_\alpha \in \mathfrak{M}} \in \underleftarrow{\lim}_{M \in \mathfrak{M}} G/M$.
		Since the bonding maps $G/M \to G/M'$ ($M \leq M'$) are the canonical projections, 
		the fact that $(g_\alpha M_\alpha)_{M_\alpha \in \mathfrak{M}}$ is in the inverse limit means that 
		\begin{equation} \label{komp} \text{ if } \ M_\alpha \geq M_\beta, \text{ then }\  g_\alpha M_\alpha \supseteq g_\beta M_\beta\end{equation}
		From that it can be easily seen that 
		the system $\{ g_\alpha M_\alpha : \ M_\alpha \in \mathfrak{M} \} \subseteq \mathcal{P}(G)$ of compact sets has the finite intersection property.

		Now, there is an element $g \in \bigcap_{M_\alpha \in \mathfrak{M}}g_\alpha M_\alpha \neq \emptyset$, then it can be easily seen that
		\[ \psi(g) = (gM_\alpha)_{M_\alpha \in \mathfrak{M}} = (g_\alpha M_\alpha)_{M_\alpha \in \mathfrak{M}}.\]
	
	The kernel of $\psi$ is $L$ indeed, since
		\[ \Ker(\psi) = \psi^{-1}((M)_{M \in \mathfrak{M}}) = \bigcap_{M \in \mathfrak{M}} M = L.\]
		 
		For the openness of $\psi$, since $\psi$ is a homomorphism between topological groups, we only have to check that for each open set $U$ containing the identity $\psi(U)$ is a neighborhood of the identity.
		Let $U$ be an open set containing the identity, and, since $\psi(U) = \psi(UL)$ (recall that $\Ker(\psi)=L$),
		 we can assume that $UL = U$ (the set $VX$ is open for arbitrary $X \subseteq G$ if $V\subseteq G$ is open).
		From $UL=U$ follows $L \subseteq U$.
		
		Now, using the compactness of the $M$-s,  we will show that at least one of the $M$-s must be contained in $U$.
		Indeed, if  the $M \setminus U$-s  are nonempty compact sets (for all $M \in \mathfrak{M}$),
		then because of the directedness, the system $\{ M \setminus U : \ M \in \mathfrak{M} \}$ of compact sets has the finite intersection property, thus
		\[ \emptyset \neq \bigcap_{M \in \mathfrak{M}} (M \setminus U) = \left( \bigcap_{M \in \mathfrak{M}} M \right) \setminus U = L \setminus U  = \emptyset. \]
		Thus there is an $M_0 \in \mathfrak{M}$ with $M_0 \subseteq U$. Then by some standard compactness argument, there is an open $U'$ with $e \in U'$
		such that 
		\begin{equation} \label{tartaos}
		U' M_0 \subseteq U. \end{equation}
		
%
		 Now we will check that 
		the basic open set 
		\[ B= \left( \underleftarrow{\lim}_{M \in \mathfrak{M}} G/M \right) \cap \left( \{ xM_0 : \ x \in U' \} \times \prod_{M \in \mathfrak{M} \setminus \{M_0\} } G/M \right) \]
		is contained in $\psi(U)$.
		Fix a $b \in B$. By the surjectivity of $\psi$,  $b= \psi(g)$ for some $g \in G$, so
		\[ \pi_{\mathfrak{M} \to \{M_0\}}(\psi(g)) = gM_0 \in \{ xM_0 : \ x \in U' \}, \]
		and then $g \in gM_0 \subseteq U'M_0 \subseteq U$ by $\eqref{tartaos}$, i.e. $g \in U$. Hence
		\[ b = \psi(g) \in \psi(U), \]
		and
		\[ \psi(U) \supseteq B,\]
		thus $\psi$ is an open mapping.

		For our second claim, observe that for any open subset $U \subseteq G/L$
		\[ \overline{\psi}(U) = \psi(\varphi_L^{-1}(U)) \]
		is open by openness of $\psi$, and the continuity of the canonical projection $\varphi_L: G \to G/L$, thus $\overline{\psi}$ is a topological group isomorphism, indeed.

\ep

\bcor \label{kov}
	Let $G$ be a locally compact FL group.
	\begin{enumerate}[(i)]
		
	\item 	\label{i}
	Then the mapping
	\[ g \mapsto (gN)_{N \in \mathfrak{N}} \]
	is an isomorphism, hence
	\[ G \simeq \underleftarrow{\lim}_{N \in \mathfrak{N}} G/N.\]
	\item \label{ii} For any directed subsystem $\mathfrak{M}$ of $\mathfrak{N}$,
	\begin{equation} \label{ize}
		\textrm{the projection	}\pi_{\mathfrak{N} \to \mathfrak{M}} \textrm{ is a surjective open mapping.}
	\end{equation}
	\end{enumerate}

\ecor
\bp

	For $\eqref{i}$ 	We only have to check the conditions of Lemma $\ref{reszrsz}$:

			We need
			$\mathfrak{N} = \{N \vartriangleleft G : \ N \textrm{ is compact,}$ $G/N$ is a Lie group with finitely many connected components $\}$
			to be directed, but $\eqref{directed}$ stated that this is true.

			 $\bigcap_{N \in \mathfrak{N}} N = \{e\}$:
			
			For any $g \in G$, if $g \neq e$, then there is an open neighborhood $U$ of the identity such that $g \notin U$. Then $G/N_U$ is a Lie group with finitely many components, i.e. $N_U \in \mathfrak{N}$.
			
Our second claim $\eqref{ii}$ easily follows from that $\pi_{\mathfrak{N} \to \mathfrak{M}}$ is the composition of the inverse of the
		 homeomorphism $\varphi : G \to \underleftarrow{\lim}_{N \in \mathfrak{N}} G/N$  ($g \mapsto (gN)_{N \in \mathfrak{N}}$),
			and the aforementioned open projection $\psi_\mathfrak{M}: G \to \underleftarrow{\lim}_{M \in \mathfrak{M}} G/M$ ($g \mapsto (gM)_{M \in \mathfrak{M}}$):
			\[ \pi_{\mathfrak{N} \to \mathfrak{M}} = \psi_\mathfrak{M} \circ \varphi^{-1}. \]

\ep

Now, we can turn to a key lemma:
\bl \label{kulcs}
	Let $(\{G_\alpha: \ \alpha \in I \}, \{ \varphi_{\alpha,\beta} : \ \alpha \preceq \beta \in I \} )$ be an inverse system of second countable topological groups.
	
	Assume that 
	a countable subset  $J_0 \subseteq I$ of the indexing poset is given, and let
	$V_k$ ($k \in \omega$) be dense open sets in $\underleftarrow{\lim}_{\alpha \in I} G_\alpha$.
	Assume that for any directed subsystem $H \subseteq I$ the projection
	\begin{equation} \label{t1}  \pi_{I \to H}: \underleftarrow{\lim}_{i \in I} G_i \to  \underleftarrow{\lim}_{i \in H} G_i \text{ is surjective.}
	\end{equation}

	Then there is a countable directed subsystem $J \subseteq I$, $J_0 \subseteq 
	J$, and a sequence  $(V_k')_{k \in \omega}$ of dense open sets in $\underleftarrow{\lim}_{\alpha \in J} G_\alpha$,
	such that for each $k \in \omega$
	\begin{equation} \label{cell} 
	  \pi_{I \to J}^{-1}(V_k') \subseteq V_k.\end{equation}
\el
\br \upshape
	Condition $\eqref{t1}$ could be easily dropped, but in our applications this condition will always hold, thus for simplicity, and for avoiding confusions  we assume it.
\er

\bp

	We can assume that $J_0 \neq \emptyset$ (by choosing an arbitrary $\alpha \in I$ to be an element of $J_0$).
	We will find a countable set $J$ such that  $J_0 \subseteq J\subseteq I$, $J$ is  directed, and closed under the following countably many operations.
	For each $\alpha \in I$ fix a countable base $\mathcal{B}_\alpha= \{B_1^\alpha, B_2^\alpha, \dots,\}$ consisting of nonempty open subsets
	  of $G_\alpha$. 
	  Recall that there is a base of $\underleftarrow{\lim}_{\alpha \in I} G_\alpha $ consisting of sets of the form 
	  \[ D_n^\alpha :=\left( B_n^\alpha \times \prod_{\beta \in I \setminus \{ \alpha \} } G_\beta \right)  \cap \underleftarrow{\lim}_{\gamma \in I} G_\gamma   \ \ (n\in \omega, \ \alpha \in I) .\]
	  (Since $B_n^\alpha$-s are nonempty open sets, by $\eqref{t1}$ the $D_n^\alpha$-s are also nonempty.)
	  Now for every $n,m \in \omega$ we define the following operation, which is a choice function.
	  \[ f_{n,m}: I \to I: \ \ 
	  \alpha \mapsto \beta \ \  \text{ s.t. } \exists k: \  D^\beta_k \subseteq  V_n \cap D_m^\alpha \]
	  (such a $\beta$ exists because each $V_n$ is a dense open set).
	  And let $g:I \times I \to I$ be another choice function such that $\alpha, \beta \preceq  g(\alpha, \beta)$, i.e. a subset 
	  of $I$ closed under 
	  the function $g$ is a directed subsystem.
	  Then, take the closure $J$ of $J_0$ under the $f_{n,m}$-s and $g$. Then, because $J$ is closed under $g$, it is a directed system.
	  
	  	  Similarly to the $D^\beta_n$-s, we define nonempty basic open sets in the limit $\underleftarrow{\lim}_{\gamma \in J} G_\gamma$ as follows.
	  \begin{equation} \label{none} E_n^\alpha :=\left( B_n^\alpha \times \prod_{\beta \in J \setminus \{ \alpha \} } G_\beta \right)  \cap \underleftarrow{\lim}_{\gamma \in J} G_\gamma  \neq \emptyset \ \ (\text{for } n\in \omega, \ \alpha \in J ), \end{equation}
	  which are nonempty because $B_n^\alpha$-s were nonempty, and using $\eqref{t1}$ with $H = \{ \alpha \}$ for each fixed $\alpha \in J$.

	  We have the
	  directed system $J \subseteq I$, and for each $n,m \in \omega$, and $\alpha \in J$, since $J$ is closed under the function $f_{n,m}$, there is a $k=k_{n,m,\alpha} \in \omega$, and 
	  $\beta= \beta_{n,m,\alpha} \in J$ such that
	  \begin{equation} \label{oo}
			 D_{k_{n,m,\alpha}}^{\beta_{n,m,\alpha}} \ \subseteq \ V_n \cap D_m^\alpha.
	  \end{equation}

	  Next we claim that for each $n,m \in \omega$, and $\alpha \in J$, if $k_{n,m,\alpha} \in \omega$, and $\beta_{n,m,\alpha} \in J$ are what $\eqref{oo}$ gives, then
	  \begin{equation} \label{ff1}
	  	\emptyset \neq E_{k_{n,m,\alpha}}^{\beta_{n,m,\alpha}} \subseteq E_m^\alpha 
	  \end{equation}   
	    and
	  \begin{equation} \label{ff2}
	    D_{k_{n,m,\alpha}}^{\beta_{n,m,\alpha}}= \pi_{I \to J }^{-1}(E_{k_{n,m,\alpha}}^{\beta_{n,m,\alpha}}) \ \  \subseteq \  \ V_n .
	    \end{equation}

		For $\eqref{ff1}$  first let $k= k_{n,m,\alpha} \in \omega$, and $\beta =\beta_{n,m,\alpha} \in J$.
		$E_{k_{n,m,\alpha}}^{\beta_{n,m,\alpha}} = E^k_\beta $ is nonempty by $\eqref{none}$.
		Reformulating
       the condition		 $D_k^\beta \subseteq D_m^\alpha$ (following from $\eqref{oo}$), we get that
		 \[ \left[ (x_\gamma)_{\gamma \in I} \in \underleftarrow{\lim}_{\gamma \in I} G_\gamma \ \wedge \ x_\beta \in B_k^\beta\right]  \ \Rightarrow \ (x_\gamma)_{\gamma \in I} \in D_m^\alpha \ \ ( \iff \  x_\alpha \in B_m^\alpha \ ).\]

		 Now notice that by our assumptions, the projection $\pi_{I \to J}: \underleftarrow{\lim}_{\gamma \in I} G_\gamma \to \underleftarrow{\lim}_{\gamma \in J} G_\gamma$ is surjective, thus for a fixed element 
		 $(x_\gamma)_{\gamma \in J}$ in the inverse limit of the smaller system there is an element $(y_\gamma)_{\gamma \in I} \in \underleftarrow{\lim}_{\gamma \in I}G_\gamma$ with $(x_\gamma)_{\gamma \in J} = \pi_{I \to J}((y_\gamma)_{\gamma \in I})$. Then
		 \[ (x_\gamma)_{\gamma \in J} \in E^\beta_k \Rightarrow x_\beta \in B_k^\beta  \Rightarrow y_\beta \in B_k^\beta \Rightarrow y_\alpha \in B_m^\alpha \Rightarrow x_\alpha \in B_m^\alpha \Rightarrow (x_\gamma)_{\gamma \in J} \in E_m^\alpha .
		  \]
		 This clearly yields $\eqref{ff1}$, as desired.
		 
		 For $\eqref{ff2}$, $D_{k_{n,m,\alpha}}^{\beta_{n,m,\alpha}} \ \  \subseteq \  \ V_n$ is true by $\eqref{oo}$, we only
		 have to check that $ D_{k_{n,m,\alpha}}^{\beta_{n,m,\alpha}}  = \pi_{I \to J }^{-1}(E_{k_{n,m,\alpha}}^{\beta_{n,m,\alpha}})$.
		 But for arbitrary $l \in \omega$ and $\delta \in J$
		 \[ D_l^\delta  = \pi_{I \to J }^{-1}(E_l^\delta), \]
		 since for 
		 $(x_\gamma)_{\gamma \in I} \in \underleftarrow{\lim}_{i \in I} G_i$ 
		 \[ (x_\gamma)_{\gamma \in I} \in D_l^\delta \iff x_\delta \in B_l^\gamma \iff
		 \]
		 \[ \iff (x_\gamma)_{\gamma \in J} = \pi_{I \to J}((x_\gamma)_{\gamma \in I}) \in E^\delta_l. \]

		 Finally, set  $V'_n = \bigcup_{m \in \omega, \alpha \in J} E^{k_{n,m,\alpha}}_{\beta_{n,m,\alpha}}$ . Since $m$ ranges over $\omega$ and $\alpha$ over $J$, using $\eqref{ff1}$ $V_n'$
		 is a dense open set. From $\eqref{ff2}$
		 \[ \pi_{I \to J}^{-1}(V'_n) = \bigcup_{m \in \omega, \alpha \in J} \pi_{I \to J}^{-1}(E_{k_{n,m,\alpha}}^{\beta_{n,m,\alpha}}) = \bigcup_{m \in \omega, \alpha \in J}  D_{k_{n,m,\alpha}}^{\beta_{n,m,\alpha}},\]
		 and using $\eqref{oo}$, 
		 \[ \bigcup_{m \in \omega, \alpha \in J}  D_{k_{n,m,\alpha}}^{\beta_{n,m,\alpha}} \subseteq V_n, \]
		 i.e. $\eqref{cell}$ holds as desired.

\ep

This lemma has the following consequence.

As before if the locally compact group $G$ is given, for any compact normal subgroup $C \vartriangleleft G$, $\varphi_C: G \to G/C$ denotes the canonical projection.

\bl \label{seged}
		Assume that  $G$ is a locally compact FL group, and let
		\[ \mathfrak{N} = \{N \vartriangleleft G : \ N \in \mathfrak{C} \ \wedge \ G/N \text{ is a Lie group with finitely many connected components} \}. \]
		
		Let a countable set $\mathfrak{J}_0= \{ N_{i}: \ i \in \omega \} \subseteq \mathfrak{N}$ be given.
		Assume that $R \subseteq G$ is co-meager. 
		Then there exist a countable set $\mathfrak{J} \supseteq \mathfrak{J}_0$, $\mathfrak{J} \subseteq \mathfrak{N}$,
		and a compact normal subgroup $K' = \bigcap_{N \in \mathfrak{J}} N$ of $G$, and a co-meager set $R' \subseteq G/K'$, such that
		 $\varphi_{K'}^{-1}(R') \subseteq R$, and $G/K'$ is an inverse limit of countable many second countable Lie groups, thus is a Polish group.
	
\el
\bp
	
	Recall that by Corollary $\ref{kov}$, $G \simeq \underleftarrow{\lim}_{N \in \mathfrak{N}} G/N$,
	where $\varphi: G \to \underleftarrow{\lim}_{N \in \mathfrak{N}} G/N$ denotes the canonical isomorphism (i.e. $g \mapsto (gN)_{N \in \mathfrak{N}}$). 
	
		Lemma $\ref{reszrsz}$ states that for each  directed subsystem $\mathfrak{M} \subseteq \mathfrak{N}$ the projection
			$\pi_{\mathfrak{N} \mapsto \mathfrak{M}}$
			is also surjective, and open.
		
	 First, we will work in $\underleftarrow{\lim}_{N \in \mathfrak{N}} G/N$.
	
	Let  $\hat{R} = \varphi(R)$ denote the corresponding co-meager subset of $\underleftarrow{\lim}_{N \in \mathfrak{N}} G/N$.
	We know that the $G/N$-s ($N \in \mathfrak{N}$) are Lie groups having finitely many connected components, thus are second countable.
	
	Let $(V_m)_{m \in \omega}$ be a sequence of dense open sets in $\underleftarrow{\lim}_{N \in \mathfrak{N}} G/N$ such that $\bigcap_{m \in \omega} V_m \subseteq \hat{R}$.
	We will apply
	Lemma $\ref{kulcs}$ with $I = \mathfrak{I}$, $J_0= \mathfrak{J}_0$ and the sequence $(V_m)_{m \in \omega}$, for what  we only have to check $\eqref{t1}$.
	Recall that, by Lemma $\ref{reszrsz}$, for each  directed subsystem $\mathfrak{M} \subseteq \mathfrak{N}$, the projection
	$\pi_{\mathfrak{N} \mapsto \mathfrak{M}}$
	is surjective.

	 Lemma $\ref{kulcs}$ gives a countable directed subset $J \subseteq I$ ($J_0 \subseteq J$), set $\mathfrak{J} =  J$.
	 Also by Lemma $\ref{kulcs}$ there is a sequence $(V_m')_{m \in \omega}$ of dense open sets in $\underleftarrow{\lim}_{N \in \mathfrak{J}} G/N$ such that
	\[ \pi_{\mathfrak{N} \to \mathfrak{J}}^{-1}(V_m') \subseteq V_m \ \ (\forall m \in \omega).\]

	Hence, if we set the co-meager set $\hat{R}': = \bigcap_{m \in \omega} V_m'$ then clearly
	\[ \pi_{\mathfrak{N} \to \mathfrak{J}}^{-1}(\hat{R}') = \pi_{\mathfrak{N} \to \mathfrak{J}}^{-1}\left(\bigcap_{m \in \omega} V_m'\right)=  \bigcap_{m \in \omega}\pi_{\mathfrak{N} \to \mathfrak{J}}^{-1}(V_m') \subseteq  \bigcap_{m \in \omega} V_m  \subseteq \hat{R} .\]

	But  $\underleftarrow{\lim}_{N \in \mathfrak{J}} G/N$ is a quotient of  $\underleftarrow{\lim}_{N \in \mathfrak{N}} G/N$, since we saw that 
	the projection $\pi_{\mathfrak{N} \to J}$ is surjective and open. Thus, by the first isomorphism theorem for topological groups (see \cite[Lemma 1.5.13.]{arkh}) $\underleftarrow{\lim}_{N \in \mathfrak{J}} G/N$ is isomorphic to
	the topological group $\left(\underleftarrow{\lim}_{N \in \mathfrak{N}} G/N\right)/ \Ker(\pi_{\mathfrak{N} \to \mathfrak{J}})$, i.e. it is
	a quotient of $\left(\underleftarrow{\lim}_{N \in \mathfrak{N}} G/N\right)$.
	Let $L:= \Ker(\pi_{\mathfrak{N} \to \mathfrak{J}})$, and
	\[ \tau_{L}: \left(\underleftarrow{\lim}_{N \in \mathfrak{N}} G/N \right) \to  \left(\underleftarrow{\lim}_{N \in \mathfrak{N}} G/N\right)/ L \]
	denote the projection, let 
	$\overline{\psi} : \left(\underleftarrow{\lim}_{N \in \mathfrak{N}} G/N\right)/ L \to \underleftarrow{\lim}_{N \in \mathfrak{J}} G/N $ denote the isomorphism for which 
	\begin{equation} \label{kompp}
	\pi_{\mathfrak{N} \to \mathfrak{J}} = \overline{\psi} \circ \tau_{L}.
	\end{equation}
	
	Now we have a 
	 co-meager set $\hat{R}'$ with $\pi_{\mathfrak{N} \to \mathfrak{J}}^{-1}(\hat{R}') \subseteq \hat{R}$ in $ \underleftarrow{\lim}_{N \in \mathfrak{J}} G/N$.

	Set $\hat{R}'' = \overline{\psi}^{-1}(\hat{R}')$, by $\eqref{kompp}$  that is   a co-meager set such that $\tau_{L}^{-1}(\hat{R}'') \subseteq \hat{R}$.
	 Now we have a normal subgroup 
	 $L=\Ker(\pi_{\mathfrak{N} \to \mathfrak{J}})$ in  $\underleftarrow{\lim}_{N \in \mathfrak{N}} G/N$, such that there is a co-meager $\hat{R}''$ in the quotient $\left(\underleftarrow{\lim}_{N \in \mathfrak{N}} G/N\right)/ L$, where 
	 \begin{equation} \label{tartalm}
		\tau_{L}^{-1}(\hat{R}'') \subseteq \hat{R} = \varphi(R).
	 \end{equation} 	
	We only have to pull pack this construction to $G$ by $\varphi$, and check that the normal subgroup $L=\Ker(\pi_{\mathfrak{N} \to \mathfrak{J}})$ is a compact normal subgroup, or equivalently its pull-back $\varphi^{-1}(\Ker(\pi_{\mathfrak{N} \to \mathfrak{J}}))$
	is compact.
	
	 Since  
	 \[ \Ker(\pi_{\mathfrak{N} \to \mathfrak{J}}) = \left( \prod_{N \in \mathfrak{J}} \{N\} \times \prod_{N \in \mathfrak{N} \setminus \mathfrak{J}} G/N \right) \cap  \underleftarrow{\lim}_{N \in \mathfrak{N}} G/N ,\]
	 we have that
	 \[ \varphi^{-1} (\Ker(\pi_{\mathfrak{N} \to \mathfrak{J}})) = \varphi^{-1}\left( \prod_{N \in \mathfrak{J}} \{N\} \times \prod_{N \in \mathfrak{N} \setminus \mathfrak{J}} G/N \right) = \bigcap_{N \in \mathfrak{J}} N ,\] 
that is indeed a compact subgroup.

Let $K' =\varphi^{-1} (\Ker(\pi_{\mathfrak{N} \to \mathfrak{J}})) = \bigcap_{N \in \mathfrak{J}} N$. 
Now, since $\varphi: G \to \underleftarrow{\lim}_{N \in \mathfrak{N}} G/N$ is an isomorphism, which maps $K'$ onto 
$L =  \Ker(\pi_{\mathfrak{N} \to \mathfrak{J}})$, there is an induced isomorphism $\overline{\varphi} : G/K' \to  \left( \underleftarrow{\lim}_{N \in \mathfrak{N}}G/N \right) / L$. Now it is easy to check that if we set $R' = \overline{\varphi}^{-1}(\hat{R}'')$, then by $\eqref{tartalm}$
\[ \varphi_{K'}(R') \subseteq R .\]


	
	Finally, it remains to check that
	\[ G/K' \simeq \left( \underleftarrow{\lim}_{N \in \mathfrak{N}} G/N \right) / \Ker(\pi_{\mathfrak{N} \to \mathfrak{J}}) \simeq \underleftarrow{\lim}_{N \in \mathfrak{J}} G/N \]
	 is an inverse limit of countably many Polish spaces.
 	 We know that
	$\underleftarrow{\lim}_{N \in \mathfrak{J}} G/N$ is an inverse limit of countably many Lie groups, each having finitely many connected components, from which each $G/N$ is a Polish group.
	Therefore $\underleftarrow{\lim}_{N \in \mathfrak{J}} G/N$ is an inverse limit of countably many Polish groups,
	i.e. it is a closed subset of the product of countably many Polish groups, thus is Polish.

\ep

\bl \label{vege}
		Assume that  $G$ is a locally compact FL group.
		Let a countable set $J_0 \subseteq \mathfrak{N}$ be given which is directed.
		Set $K: = \bigcap_{N \in J_0} N$, and
		 let $M \subseteq G/K$ be a non-meager subset (in the factor topology). Then $\varphi_{K}^{-1}(M)$ is non-meager in $G$.
\el
\bp
	On the contrary, assume that the set $R = G \setminus \varphi_{K}^{-1}(M)$ is co-meager.
	Then (by Lemma $\ref{seged}$) there is a  compact normal subgroup $K' \leq K$, and  co-meager set $R' \subseteq G/K'$ , such that 
	$G/K'$ is Polish, and
	\[ \varphi_{K'}^{-1}(R') \subseteq R ,\]
	i.e.
	\[  \varphi_{K'}^{-1}(R') \cap \varphi_{K}^{-1}(M) = \emptyset. \]
	 Then clearly (in $G/K'$)
	\begin{equation} \label{diszj}
	 R' \cap \varphi_{K/K'}^{-1}(M) = \emptyset.
	 \end{equation}

	But 
	\[ \varphi_{K/K'}: G/K' \to G/K \]
	is a surjective open mapping between Polish spaces (since by Lemma $\ref{reszrsz}$, $G/K \simeq \underleftarrow{\lim}_{N \in J_0} G/N$ is the inverse limit of countably many Polish spaces, thus is Polish), and according to \cite[Lemma 2.6 ]{levelsets} it maps a co-meager set onto a co-meager set. Thus $\varphi_{K/K'}(R') \subseteq G/K$ is co-meager. But clearly (using $\eqref{diszj}$  and $G/K'= \varphi_{K/K'}^{-1}(G/K)$)
	\[ \varphi_{K/K'}(R') \cap M = \emptyset, \]
	a contradiction.
\ep

We will need that if $H \leq G$ is an open subgroup,
$S \leq H/L$ (with some compact $L \vartriangleleft H$) is null with respect to the Haar measure of $H/L$, then $\varphi_L^{-1}(S)$ is null in $H$. (That would imply that
it is null in $G$, since for an open subgroup $H$ the restriction of the Haar measure $\mu_G|_H$ is a Haar measure of $H$.)

The following lemma can be found \cite[Sec.63, Theorem C ]{Halmoskonyv}.
\bl \label{halm1}
	Let $G$ be a locally compact group, and $\mu$ be a left-invariant Borel measure on it, which is positive on nonempty open sets, and finite on compact sets, and let $K \vartriangleleft G$ be a compact normal subgroup.
	Then the push-forward measure $B \mapsto \mu(\varphi_K^{-1}(B))$ is a left-invariant Borel measure on $G/K$, which is positive on nonempty 
	open sets, and finite on compact sets.
\el

The following lemma is well-known, but for the sake of completeness we include a proof:

\bl \label{halmos}
Let $G$ be a locally compact group, and 
\[ \mu :\mathcal{B}(G) \to [0, \infty] \]
 a left Haar Borel measure on it and let $K \vartriangleleft G$ be a compact normal subgroup, $\nu$ a left Haar Borel measure on $G/K$.
Then $\nu$ is the same (up to a multiplicative constant) as the push-forward of $\mu$, i.e. there is a constant $c>0$ such that for all $H \subseteq G/K$ Borel
\[ c\nu(H) = \mu(\varphi_K^{-1}(H)). \]

\el
\bp
	By Lemma $\ref{halm1}$, the push-forward of $\mu$ is a left-invariant Borel measure on $G/K$, which is positive on nonempty open sets, and finite on compact sets. Because the left Haar measure is unique up to a positive multiplicative constant, we only have to check the regularity of $\mu \circ \varphi_K^{-1}$.
	Let $B \subseteq G/K$ be a Borel subset, then, since taking supremum is monotonic, and continuous image of a compact set is compact,
	\[\begin{array}{rl}
	 \sup\{\mu(\varphi_K^{-1}(C))  : &  C \subseteq B \subseteq G/K , \ C \text { is  compact} \} \geq \\
    \sup\{\mu(\varphi_K^{-1}(\varphi_K(D))) : &  D \subseteq \varphi_K^{-1}(B) \subseteq G \text { is  compact} \}  = \\
     \sup\{\mu(DK)  : &  D \subseteq \varphi_K^{-1}(B) \subseteq G \text { is  compact} \} \geq \\
   \sup\{\mu(D) :&  D \subseteq \varphi_K^{-1}(B) \subseteq G \text { is  compact} \} .\end{array}
   \]
   (The first inequality holds, because if $D \subseteq \varphi_K^{-1}(B)$ is compact, then $\varphi_K(D) \subseteq B$ is compact, and the equality is by $\varphi_K^{-1}(\varphi_K(D)) = DK$.)
 And
\[  \sup\{\mu(D) : \ \   D \subseteq \varphi_K^{-1}(B) \subseteq G \text { is  compact} \}  = \mu(\varphi_K^{-1}(B)) \]
by the regularity of $\mu$.
\ep

Recall that under a left Haar measure we mean the completion of a left Haar Borel measure.
\bcor \label{nullpullback}
	Let $G$ be a locally compact group, $\mu$ is a left Haar measure on $G$. Suppose that $K \vartriangleleft G$ is a compact normal subgroup, $\nu$ is a left Haar measure on $G/K$. If $X \subseteq G/K$ is null wrt. $\nu$, then $\varphi_K^{-1}(X)$ is null wrt. $\mu$.
\ecor
\bp
Let $B \supseteq X$ be a null Borel set in $G/K$. Then, applying Lemma $\ref{halmos}$ for the Borel measures $\mu|_{\mathcal{B}(G)}$ and $\nu|_{\mathcal{B}(G/K)}$ we obtain that (by the continuity of $\varphi_K$)  
$\mu|_{\mathcal{B}(G)}(\varphi_K^{-1}(B)) = 0$, hence
\[ \mu(\varphi_K^{-1}(X)) \leq \mu(\varphi_K^{-1}(B)) = 0 .\]
\ep

Before the construction for profinite groups, and connected Lie groups,  we will need the following technical lemma about co-meager sets in compact metric spaces (the lemma uses the idea of the well known characterization of co-meager sets in the Cantor space, see \cite[Thm. 2.2.4 ]{Bart}).

Before stating the lemma, recall that in a metric space $(X,d)$, the diameter of a subset $S \subseteq X$ is defined as follows
\[ \diam(S) =  \sup\{ d(s,s'): \ s,s' \in S \} .\]
\bl \label{rezid}
Assume that $(X,d)$ is a compact metric space, and let $(M_i)_{i \in \omega}$ be a sequence of finite sets, let $R \subseteq X$ be a co-meager set.

Assume that  the system of compact sets $(C_{p_0p_1 \dots p_i})_{p_0p_1p_2\dots p_i \in \prod_{j=0}^{i}M_j},  (i \in \omega )$ fulfills the following conditions:

\begin{enumerate}
	\item \label{atm}	$|\bigcap_{i \in \omega} C_{p_0p_1 \dots p_i }| = 1$  for each  $p_0 p_{1} \dots p_j  \dots \in \prod_{j \in \omega} M_j  $,
	
	\item \label{resz} $C_{p_0p_1 \dots p_n p_{n+1} } \subseteq C_{p_0p_1 \dots p_n}$ for each $n \geq 0$ and $C_{p_0p_1\dots p_n p_{n+1}}$,
	
	\item  $ \inte(C_{p_0p_1\dots p_i}) \neq \emptyset$ for each $C_{p_0p_1\dots p_i}$,
	
	\item \label{un}
	for each $C_{p_0p_1\dots p_i}$
         \[ C_{p_0 p_1\dots p_i} = \bigcup_{j \in M_{i+1}} \overline{\inte(C_{p_0 p_1\dots p_i j})}  .\]
\end{enumerate}
Then there is an increasing sequence $0=n_0 < n_1 <  \dots n_i < \dots$ of integers, and an infinite  sequence 
\[ r:=(r_j)_{j \in \omega} \in \prod_{j \in \omega}M_j \]
such that for each sequence $s=(s_j)_{j \in \omega} \in \prod_{j \in \omega} M_j$, if $ \{ i \in \omega : s|_{[n_i,n_{i+1})} = r|_{[n_i,n_{i+1})]} \}$ is infinite, then $\bigcap_{j \in \omega} C_{s_0s_1\dots s_j}   \subseteq R$.
\el
\bp
     Let $\bigcap_{k \in \omega} U_k$ be the intersection of dense open sets such that $\bigcap_{k \in \omega} U_k \subseteq R$.
      By replacing $U_k$ with $U_0 \cap U_1 \cap \dots \cap U_k$ 
	we can assume that $(U_k)_{k \in \omega}$ is a decreasing sequence. We construct the sequences $(n_i)_{i \in \omega}$  and $(r_i)_{i \in \omega}$ simultaneously by induction on $i$ as follows.
	Before the $i+1$-th step the initial segments $0= n_0<n_1 <n_2 < \dots <n_{i}$, and $r_0r_1\dots r_{n_{i}-1}$ are already defined, and we will choose $n_{i+1}$, and $r_{n_{i}}r_{n_{i}+1}\dots r_{n_{i+1}-1}$ such that
		
	\begin{equation}
		 \label{ffff}
		 \begin{array}{l}
		 C_{p_0p_1 \dots p_{n_i-1}r_{n_i}r_{n_i+1}r_{n_i+2}\dots r_{n_{i+1}-1}} \subseteq U_{i+1} \\
		 \text{ for each } p_0p_1\dots p_{n_i-1} \in \prod_{j <n_i} M_j.
		 		 \end{array}
	\end{equation} 

	First, we show that we can choose $n_{i+1}$  and $r_{n_{i}}r_{n_{i}+1}\dots r_{n_{i+1}-1}$ so that $\eqref{ffff}$ holds.

	Assume that $n_0, n_1,n_2, \dots ,n_i$ are defined. Now enumerate the set
	\[ \{ C_{p_0p_1 \dots p_{n_i-1}}: \ \ p_0p_1 \dots p_{n_i-1} \in \prod_{j<n_i} M_j  \} = \{ D_1,D_2, \dots, D_{|M_0|\cdot|M_1|\cdot|M_2| \cdot \dots \cdot |M_{n_i-1}|  } \}. \]
	For each $D_l$ define the sequences  $\underline{r}^l= r^l_1r^l_2\dots r^l_{h_l}$
	 by induction on $l$, such that if $D_l = C_{p_0p_1 \dots p_{n_i-1}}$, then for the concatenation $\underline{r}^1 \underline{r}^2 \dots \underline{r}^l$ 
	of the $\underline{r}^j$-s ($j \leq l$)
	\begin{equation} \label{kozt}
	 C_{p_0p_1 \dots p_{n_i-1}\underline{r}^1 \underline{r}^2 \dots \underline{r}^l} \subseteq U_{i+1} .
	 \end{equation}
	Before constructing the $\underline{r}^j$-s ($j \leq |M_0|\cdot |M_1|\cdot |M_2| \cdot \dots \cdot |M_{n_i-1}|$), we prove that if
 we could manage to choose such $\underline{r}^l$-s, then 
	setting 
	\begin{itemize}
		\item $n_{i+1} = n_i + \length(\underline{r}^1 \underline{r}^2 \dots \underline{r}^{|M_0|\cdot |M_1|\cdot |M_2| \cdot \dots \cdot |M_{n_i-1}|} )$ , and
		\item $r_{n_i}r_{n_i+1} \dots r_{n_{i+1}-1} = \underline{r}^1 \underline{r}^2 \dots \underline{r}^{|M_0|\cdot |M_1|\cdot |M_2| \cdot \dots \cdot  |M_{n_i-1}|}$
	\end{itemize}
		would yield $\eqref{ffff}$.
		 Since for $D_l = C_{p_0p_1 \dots p_{n_i-1}}$ $\eqref{resz}$ in the conditions of the lemma and $\eqref{kozt}$ implies that
		\[ C_{p_0p_1 \dots p_{n_i-1}\underline{r}^1 \underline{r}^2 \dots \underline{r}^{|M_0|\cdot |M_1|\cdot |M_2|\cdot \dots \cdot |M_{n_i-1}|}} \subseteq 
		C_{p_0p_1 \dots p_{n_i-1}\underline{r}^1 \underline{r}^2 \dots \underline{r}^l} \subseteq U_{i+1} , \] 
		i.e. $\eqref{ffff}$ holds indeed.
	
	Turning to the construction of the $\underline{r}^l$-s, assume that $\underline{r}^1, \underline{r}^2, \dots \underline{r}^{l-1}$ are defined. If  $D_l = C_{p_0p_1 \dots p_{n_i-1}}$, then, because
	\[ \inte( C_{p_0p_1 \dots p_{n_i-1} \underline{r}^1 \underline{r}^2 \dots \underline{r}^{l-1}}) \neq \emptyset \]
	and $U_{i+1}$ is a dense open set,
	 there is an $x \in X$, and an open ball 
	 \[ B(x,r) = \{ y \in X: d(x,y)<r \} \subseteq  C_{p_0p_1 \dots p_{n_i-1} \underline{r}^1 \underline{r}^2 \dots \underline{r}^{l-1}} \cap U_{i+1} . \]
	 Now we can use $\eqref{un}$ from our assumptions, thus there is
	 an $r^l_0$ such that for the open set $V_0 = \inte(C_{p_0p_1 \dots p_{n_i-1} \underline{r}^1 \underline{r}^2 \dots \underline{r}^{l-1}r^l_0})$
	 \[ V_0 \cap B\left(x,\frac{r}{2}\right) \neq \emptyset .\]
	 Similarly, by $\eqref{un}$ of our conditions, there are infinite sequences $r^l_0, r^l_1, r^l_2, \dots , r^l_k, \dots$  and 
	 $V_0, V_1, \dots $, where $V_k = \inte(C_{p_0p_1 \dots p_{n_i-1} \underline{r}^1 \underline{r}^2 \dots \underline{r}^{l-1} r^l_0r^l_1\dots r^l_k})$ such that
	  	 \[ V_k \cap  B(x,\frac{r}{2}) \neq \emptyset .\]
	  And using $\eqref{atm}$ (from the lemma) we know that
	 \[  \lim_{k \to \infty}  \diam(C_{p_0p_1 \dots p_{n_i-1} \underline{r}^1 \underline{r}^2 \dots \underline{r}^{l-1} r^l_0r^l_1\dots r^l_k} )= 0, \]
	  thus there is an index $k_0$ such that 
	  \[ \diam(V_{k_0}) \leq \diam(C_{p_0p_1 \dots p_{n_i-1} \underline{r}^1 \underline{r}^2 \dots \underline{r}^{l-1} r^l_0r^l_1\dots r^l_{k_0}}) \leq  \frac{r}{4} ,\]
	  from which 
	  \[ C_{p_0p_1 \dots p_{n_i-1} \underline{r}^1 \underline{r}^2 \dots \underline{r}^{l-1}r^l_0r^l_1\dots r^l_{k_0}} \subseteq \overline{B(x,\frac{3r}{4})}  ,\]
	  thus
	  \[ C_{p_0p_1 \dots p_{n_i-1} \underline{r}^1 \underline{r}^2 \dots \underline{r}^{l-1} r^l_0r^l_1\dots r^l_{k_0}} \subseteq B(x,r) \subseteq U_{i+1}. \]
	  This finishes the construction of $(n_i)_{i \in \omega}$ and $(r_i)_{i \in \omega}$.
	  
	  Now, as the $r_i$-s, and $n_i$-s are already constructed, we only have to check that for any sequence
	  $(s_j)_{j \in \omega} \in \prod_{j \in \omega} M_j$ if $ \{ j \in \omega : r_{|[n_j,n_{j+1})} = s_{|[n_j,n_{j+1})]} \}$ is infinite, then $\bigcap_{i \in \omega} C_{s_0s_1\dots s_i}   \subseteq R$. Fix such a sequence $s = (s_j)_{j \in \omega}$.
	  
	  If $k \in \{ j \in \omega : r_{|[n_j,n_{j+1})} = s_{|[n_j,n_{j+1})]} \}$, then, using $\eqref{ffff}$:
	  \[ C_{s_0s_1 \dots s_{n_k-1}r_{n_k}r_{n_k+1}r_{n_k+2}\dots r_{n_{k+1}-1}} \subseteq U_{k+1} . \]
	  Thus, if $x \in \bigcap_{i \in \omega} C_{s_0s_1\dots s_i}$, $k \in \{ j \in \omega : r_{|[n_j,n_{j+1})} = s_{|[n_j,n_{j+1})]} \}$ are fixed, then
	   \[ x  \in C_{s_0s_1 \dots s_{n_{k+1}-1}} =  C_{s_0s_1 \dots s_{n_k-1}r_{n_k}r_{n_k+1}r_{n_k+2}\dots r_{n_{k+1}-1}} \subseteq U_{k+1} ,\]
	   hence $x \in U_{k+1}$ for infinitely many $k$-s. But recall that the $U_k$-s form a decreasing sequence, thus
	   \[ x \in \bigcap_{k \in \omega} U_k \subseteq R .\]

\ep

\bprop \label{profinite}
	Let $G$ be the inverse limit of the following inverse system, consisting of finite groups.
	\[ \{ G_i , \varphi_{i,j}:G_j \to G_i, \ i\leq j \in \omega \}  ,\]
	i.e. 
	\[ G = \{  (g_i)_{i \in \omega} : \ \forall j \ \varphi_{j,j+1}(g_{j+1}) = g_j \} \leq \prod_{j \in \omega} G_j, \]
	where each $|G_{j+1}|/|G_j|>1$. Then there exists a subgroup $S \leq G$ that is null but non-meager.
\eprop 
\bp
	
	Let $\varphi_i: G \to G_i$ denote the canonical projection.
	Since, by  \cite[Lemma 1.1.9 ]{profinite}, for an arbitrary increasing sequence $n_1 < n_2 < \dots $ of positive integers
	\begin{equation} \label{thinning} \underleftarrow{\lim}_{i \in \omega} G_i \simeq \underleftarrow{\lim}_{i \in \omega} G_{n_i}, \end{equation}
	we can assume that  the sequence 
	\[ \left(m_i = \frac{|G_i|}{|G_{i-1}|}\right)_{1 \leq i \in \omega} \]
	 can grow as rapidly as we want. Later we will recursively thin out the sequence of $G_i$-s. Another corollary of this lemma is that we can assume that $G_0 = \{e \}$ is the trivial group. Define $m_0 = 1$, thus from now $\dom(m)= \omega$.
	
	For each $g \in G_{i-1}$ ($i \geq 1$) fix an eumeration of the set
	\[ \varphi_{i-1,i}^{-1}(g) = \{ g^{(1)},g^{(2)}, \dots, g^{(m_i)} \}, \]
	so that 
	\begin{equation} \label{egyseg}  e_{G_{i-1}}^{(1)} = e_{G_{i}}  \end{equation}
	holds.
	For $i \geq 1$ let $G_i^{(j)} \subseteq G_i$ denote the set $\{ g^{(j)} : g \in G_{i-1} \}$, i.e. for each $g \in G_{i-1}$ it contains exactly one element from $\varphi_{i-1,i}^{-1}(g)$, we have a partition of $G_i = \cup_{j=1}^{m_i} G_i^{(j)}$. For the trivial group $G_0$ define $G_0^{(1)} = \{e_{G_0}\}$.
	
	Fixing these enumerations, each element of $G$ can be uniquely identified with an element of $\prod_{i \in \omega} \{1,2, \dots,m_i \}$ in the following way.
	Let $\psi_i : G_i \to \{ 1,2, \dots, m_i\}$ denote the mapping 
	\begin{equation} \label{psidef}
	  g \mapsto k \text{ iff }  g \in G_i^{(k)}, \end{equation} 
	  i.e. for $i \geq 1$ $g = (\varphi_{i-1,i}(g))^{(k)}$. Then consider the product $\prod_{i \in \omega} \{1,2, \dots,m_i \}$ with the product topology, and the mapping 
	\[ \psi: \underleftarrow{\lim}_{i \in \omega} G_i \to \prod_{i \in \omega} \{1,2, \dots,m_i \} \]
	\[ (g_i)_{i \in \omega} \mapsto (\psi_i(g_i))_{i \in \omega} .\]
	Notice that for $e_G = (e_{G_i})_{i \in \omega} \in G$  $\eqref{egyseg}$ implies
	\begin{equation} \label{ee}  \psi(e_G) = (1)_{i \in \omega} \end{equation}
	\bc
		The mapping $\psi: G \to  \prod_{j \in \omega} \{1,2,\dots,m_j\} $ is a  homeomorphism.
		Furthermore, considering $G$ with the probability left Haar measure $\mu_G$ (i.e. $\mu_G(G) = 1$), and the product space with the canonical measure $\nu$, i.e. \[ \nu(\{k_0\} \times \{k_1 \} \times \dots \{ k_i\} \times \prod_{j > i} \{1,2,\dots,m_j\}) = \frac{1}{m_0} \cdot \frac{1}{m_1} \cdot \dots \cdot \frac{1}{m_i} ,\]
		$\psi$ is measure preserving. 
		
	\ec
	\bp
		It is straightforward to check that $\psi$ is a bijection.
		
		Fix an element $g \in G_i$ of the $i$-th group and
		the corresponding basic set
		\[ B = \left(\{g\} \times \prod_{j \in \omega, j \neq i} G_j \right) \cap \underleftarrow{\lim}_{j \in \omega} G_j \]
		 in $\underleftarrow{\lim}_{j \in \omega} G_j$.
		
		Then consider the following sequence
		\[ h_i=g, \  h_{i-1}= \varphi_{i-1,i}(g), \dots, \ h_0 = \varphi_{0,i}(g) \]
		from which we can define $k_j$-s ($j \leq i$, $1 \leq k_j \leq m_j$) such that $h_j = h_{j-1}^{(k_j)}$. Then one can easily see that 
		\[ \psi(B) = \{k_0\} \times \{k_1 \} \times \dots \times \{ k_i\} \times \prod_{j > i} \{1,2,\dots,m_j\} .\]
		As $\psi$ is an open bijection onto a compact set, it is a homeomorphism.

		It is left to show that $\psi$ is measure preserving. Recall that $\varphi_i : G= \underleftarrow{\lim}_{j \in \omega} G_j \to G_i$ is the canonical projection. For the basic open set
		\[ B= \{ l_0 \} \times \{ l_1 \} \times \{l_2\} \times \dots \times \{ l_i \} \times \prod_{j>i} \{1,2,\dots,m_j\}, \]
		let $g_i \in G_i$ be such that 
		\[ \psi^{-1}(B)=  \left( \{g_i\} \times \prod_{j \in \omega, \ j \neq i} G_j  \right) \cap \underleftarrow{\lim}_{j \in \omega} G_j = \varphi_i^{-1}(g_i)   . \]
		Now, observe that if $g_i,h_i \in G_i$, and $h \in G$ is such that $\varphi_i(h)=h_i$, then $h \varphi^{-1}_i(g_i) =\varphi_i^{-1}(g_ih_i)$, 
		i.e. sets of the form $\varphi_i^{-1}(x)$ ($x \in G_i$) are translates of each other. Using this, by the left-invariance of $\mu_G$
		\[ \mu_G(\varphi_i^{-1}(g_i)) = \frac{1}{|G_i|}   =\frac{1}{m_0} \cdot \frac{1}{m_1} \cdot  \dots  \cdot \frac{1}{m_i}  = \nu(B). \]
		Thus, we got that $\nu$ and $\mu_G \circ \psi^{-1}$ coincide on sets of the form $\{k_0\} \times \{k_1 \} \times \dots \{ k_i\} \times \prod_{j > i} \{1,2,\dots,m_j\}$ (i.e. on basic open sets), from which they coincide on the generated ring, since every finite union can be written as a disjoint union. Now, from \cite[Sec. 13, Thm A ]{Halmoskonyv}, $\nu$ and $\mu_G \circ \psi^{-1}$ coincide on the generated $\sigma$-ring, that is
		the $\sigma$-algebra of the Borel sets, since the whole product $ \prod_{j \in \omega} \{1,2,\dots,m_j\}$ was contained in the ring.
	\ep

	For an arbitrary set $X$, let $\mathcal{P}(X)$ denote the power set of $X$.
	Now fix an $i \in \omega$.  We define the following operations:
	\[ \mathcal{F}_i: \mathcal{P}(G_i) \to \mathcal{P}(G_i) : 
	\ \  H \mapsto  HH \cup (HH)^{-1} \]
	\[ \mathcal{G}_i: \mathcal{P}(G_i) \to \mathcal{P}(G_i) : \ \ \{ g_1^{(j_1)},g_2^{(j_2)}, \dots, g_n^{(j_n)} \} \mapsto  \{ g^{(j_k)} : \ g \in G_{i-1}, \ k\leq n \} .\]
 So $\mathcal{G}_i(H)$ is the minimal cover of $H$ by the union of sets of the form $G_i^{(j)}$.

  It is obvious that for $H \subseteq G_i$ we have
  \[ |\mathcal{F}_i(H)| \leq 2 |H|^2 \]
  and
  \[ |\mathcal{G}_i(H)| \leq |H| |G_{i-1}|, \]
  thus
 \begin{equation} \label{upb} |\mathcal{G}_i\circ\mathcal{F}_i(H)| \leq 2|H|^2 |G_{i-1}|. \end{equation}
 From these inequalities it can be easily seen that for each $H \subseteq G_i$ there is an upper bound for  $|(\mathcal{G}_i\circ\mathcal{F}_i)^j(H)|$ depending only on $|H|$, $|G_{i-1}|$, and $j$.
 
  Now, for all $i,j$ we will define sets $B_i^j \subseteq \{ 1,2,  \dots m_i \}$  such that
  \begin{equation} \label{Bszorz}
  \text{if } k,k' \in B_i^j \ \Rightarrow \ \ G_i^{(k)} \cdot G_i^{(k')} \ \subseteq \cup \{ G_i^{(l)}  : \ l \in B_i^{j+1} \}
  \end{equation}
  and
  \begin{equation} \label{Binv}
  \text{if } k \in B_i^j \ \Rightarrow \ \ \left(G_i^{(k)}\right)^{-1} \ \subseteq \cup \{ G_i^{(l)}  : \ l \in B_i^{j+1} \}.
  \end{equation}

 Let $A_i^0 = \mathcal{G}_i(\{e_{G_i}\} ) \subseteq G_i$, and define  $A_i^j = (\mathcal{G}_i \circ \mathcal{F}_i)^j(A_0)$ for each $j \leq i$, and let $A_i^{i+1} = G_i$. 
 For every $j$, as $A_i^j$ is the union of sets of the form $G_i^{(k)}$, let 
 \begin{equation} \label{BAof} B_i^j = \{ k : \ k \leq m_i, G_i^{(k)}\subseteq A^j_i \} \subseteq \{ 1,2, \dots, m_i \} \end{equation}
	($B^i_j$ stores the information: which elements of the partition $G_i = G_i^{(1)} \cup G_i^{(2)} \cup \dots \cup G_i^{(m_i)}$ are contained in $A_i^j$, and $B_i^j=\{1,2,\dots,m_i\}$ holds for each $j > i$). 
	By the construction of the $A_i^j$-s, 
	\begin{equation} \label{Aszorz}
	 A_i^j A_i^j \subseteq A_i^{j+1} ,
	\end{equation}
	and
	\begin{equation} \label{Ainv}
	(A_i^j)^{-1} \subseteq A_i^{j+1} .
	\end{equation}
	Moreover, this and $e_{G_i} \in A_i^0$ implies that for each $j$ we have $e_{G_i} \in A_i^j$, hence  
	\begin{equation} \label{Amon} A_i^j \subseteq A_i^{j+1}  \ (j \in \omega). \end{equation}
	Thus it is straightforward to check that $\eqref{BAof}$ implies that $B_i^j$-s  satisfy  $\eqref{Bszorz}$ and $\eqref{Binv}$, and
	\begin{equation}
	 \label{Bmon} B_i^j \subseteq B_i^{j+1}  \ (j \in \omega).
	\end{equation}
	
	Next, we prove that by thinning out the sequence of the groups, we can assume that 
	\begin{equation} \label{korl} |B_i^i| i^2 \leq m_i .\end{equation}
	(Recall that by $\eqref{thinning}$, thinning out the sequence of the groups yields the same topological group.)
	\bc
		If the $\omega$-type inverse system of finite groups $H_0 = \{e \}$, $H_1$, $H_2$, $\dots$ is given (where each $H_i$ is a surjective homomorphic image of $H_{i+1}$, $|H_{i+1}|/|H_i| >1$), then one can construct a subsequence
		$(G_i = H_{n_i})_{i \in \omega}$ (where $n_0 = 0$), such that $|B_i^i| i^2 \leq m_i$ holds with the $B_i^j$-s and $m_i$-s defined for the system $(G_i)_{i \in \omega}$. 
	\ec
\bp The recursion goes as follows. Let $n_0 = 0$, and if $n_0$, $n_1$, $\dots$, $n_{i-1}$ are defined,
	then recall that for any $j \leq i$  there is an upper bound 
	\begin{equation} \label{fkorl}
	 d \geq |A_i^j| = |(\mathcal{G}_i\circ\mathcal{F}_i)^j(\{e\})|,
	\end{equation}
	 where $d$ depends only on $|G_{i-1}| = |H_{n_{i-1}}|$ and $j$, no matter what $G_i$ is (by $\eqref{upb}$). Thus we can choose $n_i$ so that 
	 \begin{equation} \label{dkorl} d i^2 \leq \frac{|H_{n_i}|}{|H_{n_{i-1}}|} \end{equation}
	 holds.
	Then, if $G_i = H_{n_i}$, using $\eqref{dkorl}$, $\eqref{fkorl}$ and the equality $|A_i^i| = |B_i^i| |G_{i-1}|$
	\[ m_i = \frac{|G_i|}{|G_{i-1}|}  = \frac{|H_{n_i}|}{|H_{n_{i-1}}|}  \geq di^2 \geq 	|A_i^i| i^2  \geq  |B_i^i| i^2 \]
   thus $\eqref{korl}$ holds, indeed.
	\ep

	Now we can construct the desired subgroup. Fix a  non-principal ultrafilter $\mathcal{U}$ on $\omega$.
	Define the subset $S' \subseteq \prod_{i \in \omega} \{ 1,2,\dots, m_i\}$ as follows.
	For each $t \in \prod_{i \in \omega} \{1,2, \dots,m_i \} $
	\[ t \in S' \ \iff \ \exists n : \ \{ i \in \omega: \  t_i \in B_i^n \} \in \mathcal{U}. \]
	Let $S = \psi^{-1}(S')$.
	We have to check that
	\begin{enumerate}[(i)]
		\item $S$ is a subgroup.
	\end{enumerate}	
		For $g = (g_i)_{i \in \omega} \in S$, let $s = \psi(g) \in \prod_{i \in \omega} \{1,2,\dots, m_i \}$ (which means that 
		\[ g = (g_i)_{i \in \omega} \ \ g_i^{(s_{i+1})}=g_{i+1},\  \forall i \in \omega \]
		by $\eqref{psidef}$, where we defined the $\psi_i$-s).
		 Then obviously
		 \[ g^{-1} =  (g_i^{-1})_{i \in \omega} , \]
		 and letting $ s'=\psi(g^{-1})$, again from the definition of $\psi_i$-s
		 \begin{equation} \label{ginverze}
			(g_{i-1}^{-1})^{(s'_i)} = g_i^{-1}.
		 \end{equation}
		Choose a positive integer $n_s$ such that $U_s = \{ i \in \omega: \ s_i \in B_i^{n_s} \} \in \mathcal{U}$. 
		Now, for $i \in \omega$ (using $g_{i-1}^{(s_i)} = g_i$ together with $\eqref{BAof}$, then $\eqref{Aszorz}$, finally $\eqref{ginverze}$  and $\eqref{BAof}$)
		\[ i \in U_s  \iff \ s_i \in B_i^{n_s}  \iff \ g_i \in A_i^{n_s}  \Rightarrow \ g_i^{-1} \in A_i^{n_s+1}  \iff \ s'_i \in B_i^{n_s+1}, \]
		which means that 
		 $n_s+1$ and $U_s$ witness that $s' \in S'$, i.e. 
		 $g^{-1} \in S$.
		 
		In order to show that $S$ is closed under multiplication, let $h \in S$ be another element, and $t= \psi(h)$
		(i.e. $h = (h_i)_{i \in \omega}$, where $h_i \in G_i$, $h_i^{(t_{i+1})}=h_{i+1}$),
		  let $n_t$ denote a positive integer for which $U_t = \{ i : \  t_i \in B_i^{n_t} \} \in \mathcal{U}$,
		and let $n = \max\{n_s,n_t\}$, $U = U_s \cap U_t$.
		Suppose that $gh = (g_ih_i)_{i \in \omega}$ is represented by the sequence of integers
		$\psi(gh) = (v_i)_{i \in \omega}$. Then  by $\eqref{Amon}$
		 \[ A^{n_s}_i, A^{n_t}_i \subseteq A^n_i, \]
		 thus for any $i \in U$
		  \[  g_i,h_i \in A^n_i .\]
		  Using $\eqref{Aszorz}$ we conclude that for $i \in U$ $ g_ih_i \in A^{n+1}_i$ holds,  
		  and hence $v_i \in B_i^{n+1}$.
	    Thus $n+1$ and $U $ witness that $gh \in S$.
	 \begin{enumerate}[(ii)]   
		\item $S$ is null.
	 \end{enumerate}	
		Since the fixed bijection $\psi$ between $G$ and $\prod_{i \in \omega} \{1,2, \dots,m_i\}$ is a measure preserving homeomorphism,
			we can work in this product space, and it is enough to prove that $S'$ is null.
		\[ S' \subseteq \bigcup_{n \in \omega} \{ s : \ \{ k \in \omega : s_k \in B^n_k \} \text{ is infinite} \} = \bigcup_{n \in \omega} \bigcap_{j\in \omega} \bigcup_{k \geq j} \{s: \ s_k \in B_k^n  \}.  \]
		By $\sigma$-additivity, showing that for abitrary fixed $n$ the set $ \bigcap_{j \in \omega} \bigcup_{k \geq j} \{s: \ s_k \in B_k^n  \}$ is null  will suffice.
		
		Recall that if $l \leq l'$, then (by $\eqref{Bmon}$)  $B^l_k \subseteq B^{l'}_k$,
		and by $\eqref{korl}$
		 $|B_k^k |k^2 \leq m_k$, hence $\nu(\{s: \ s_k \in B^k_k\}) = \frac{|B_k^k|}{m_k} \leq \frac{1}{k^2}$.
		 So if $k \geq n$, then
		\[  \nu(\{s: \ s_k \in B^n_k  \}) \leq \frac{1}{k^2}. \]
		Then, for arbitrary $n,i$ if $i\geq n$,
		 \[ \nu\left(\bigcap_{j \in \omega} \bigcup_{k \geq j} \{s: \ s_k \in B^n_k  \}\right) \leq \nu\left( \bigcup_{k \geq i} \{s: \ s_k \in B^n_k  \}\right) \leq \sum_{k=i}^\infty \frac{1}{k^2} .\]
		 Thus letting $i$ tend to infinity we obtain 
		 \[ \nu\left(\bigcap_{j \in \omega} \bigcup_{k \geq j} \{s: \ s_k \in B^n_k  \}\right) =0 .\]
		 
	\begin{enumerate}[(iii)]	
		\item $S$ is not meager.
	\end{enumerate}
			We work in the product space $\prod_{i \in \omega} \{1,2,\dots,m_i \}$  again, i.e. we prove that $S'$ is of second category.
			Let  $R$ be a co-meager set in $\prod_{i \in \omega} \{1,2,\dots,m_i \}$, we will prove that $S' \cap R$ is nonempty.
			
			Since  $X=\prod_{i \in \omega} \{1,2,\dots,m_i \}$
			is a Polish space, thus there is a compatible complete metric $d$.
			The sets of the form
			\[ C_{p_0p_1p_2\dots p_i} = \{p_0\} \times \{p_1 \} \times \dots \times \{p_i\} \times \prod_{j \in \omega, j > i} \{1,2,\dots, m_j\} \]
			are compact open sets, and for each infinite sequence $(p_i)_{i \in \omega} $
			\[ \bigcap_{j \in \omega} C_{p_0p_1p_2\dots p_j} = \{ (p_i)_{i \in \omega}\} ,\]
			therefore we can apply Lemma $\ref{rezid}$ for $M_j = \{1,2,3 \dots, m_j \}$ ($j \in \omega$): There is an increasing sequence $0= n_0< n_1 <n_2 < \dots$, and an $r \in \prod_{i \in \omega} \{1,2,\dots,m_i \}$ such that for any $s \in \prod_{i \in \omega} \{1,2,\dots,m_i \}$ if 
			\[ \{ j : s|_{[n_j,n_{j+1})} = r|_{[n_j,n_{j+1})} \} \textrm{ is infinite }, \] 
			then 
			\[ \bigcap_{i \in \omega} C_{s_0s_1 \dots s_i} =    \bigcap_{i \in \omega}\left( \{s_0\} \times \{s_1\} \times \dots \times \{s_i\} \times \prod_{j > i }\{1,2,\dots,m_j\} \right) \subseteq R  \]
			i.e. $s = (s_i)_{i \in \omega} \in R$.
			
			Now one of the following sets:
			\[ \bigcup_{j \in\omega} [n_{2j},n_{2j+1}) , \ \ \bigcup_{j \in\omega} [n_{2j+1},n_{2j+2}) \]		
			is in $\mathcal{U}$, let $U$ denote that set. Then the following sequence
			\[ s_i = \begin{cases} 
			1 & \mbox{if } i \in U \\
			r_i & \mbox{if } i \notin U  \end{cases} \]
			is obviously in $R$, because there are infinitely many intervals $[n_j,n_{j+1})$ contained in $\omega \setminus U$.
			 On the other hand, by  $\eqref{ee}$ we know that 
			$\psi(e)$ is the constant $1$ sequence in $\prod_{i \in \omega} \{1,2,\dots,m_i \}$.
			 Therefore $B^0_i= \{ 1 \}$ for all $i$, thus
			 \[ \{ i: \ s_i \in B_i^0 \} \supseteq U \in  \mathcal{U}, \]
			 from which $s \in S'$. We have that $S'$ is not disjoint from $R$.
\ep

Now we can turn to the case of Lie groups.

\bprop \label{Lieeset}
	Let $G$ be a second countable Lie group of positive dimension. Then there is a subgroup $S \leq G$ that is null, but non-meager.
\eprop
\bp
	Let $d > 0$ denote the dimension of $G$. From now on we fix a left Haar measure $\mu$  of $G$. $ [0,1]^d$ denotes  the $d$-dimensional unit cube.
	Second countable Lie groups are Lindelöf spaces, thus there are compact sets $Q_i \subseteq G$ ($i \in \omega$) homeomorphic to $[0,1]^d$ (by  homeomorphisms $\varphi_i: Q_i \to [0,1]^d$), such that $\cup_{i \in \omega} \varphi_i^{-1}((0,1)^d)$ covers $G$. Then also the $Q_i = \dom(\varphi_i) = \varphi_i^{-1}([0,1]^d)$-s (which may be overlapping) cover the group.
	We can assume that 
	\begin{equation} 
		\inte(Q_i) = \inte(\dom(\varphi_i)) = \varphi_i^{-1}((0,1)^d), \end{equation}
		in particular
		\begin{equation} \label{inpart}
		   \inte(\dom(\varphi_i)) \text{ is dense in }\dom(\varphi_i) = Q_i
		\end{equation}
			(because $(0,1)^d$ is dense in $[0,1]^d$ and $\varphi_i$ is a homeomorphism),
		and
	\begin{equation} \label{egysegegyben}
			e \in \varphi_0^{-1}((0,1)^d).
	\end{equation}
	Let $\mathfrak{C}^{i} = \{ C^i_1, C^i_2, \dots, C^i_{2^{di}}\} $ denote the $i$-th generation dyadic subdivision of the cube $[0,1]^d$, i.e. each $C^i_j$ is of the form $[\frac{k_1}{2^i}, \frac{k_1+1}{2^i}] \times [\frac{k_2}{2^i},\frac{k_2+1}{2^i}] \times \dots \times[\frac{k_d}{2^i},\frac{k_d+1}{2^i}]$, where each $k_l \in \{ 0,1, \dots 2^i-1\}  $.
	Recall that for a given set $X \subseteq G$ $\varphi_k[X]$ denotes the image of $X$ under $\varphi_k$, i.e.
	\[ \varphi_k[X] = \{ \varphi_k(x): \ x \in X \cap \dom(\varphi_k) \} .\]
	
	We will define a strictly increasing sequence $(m_i)_{i \in \omega}$, and for each $i$ a sequence $(\mathcal{D}^i_j)_{j \in \omega}$ of subsets of $\mathfrak{C}^{m_i}$,
	such that for every $i \in \omega$
	\begin{enumerate}[(i)]
		\item \label{Dmon} $\cup \mathcal{D}^i_j \subseteq \cup \mathcal{D}^i_{j+1}$ for each $j \in \omega$ (or equivalently $(\mathcal{D}^i_j)_{j \in \omega}$ is increasing),
		\item \label{De} $e \in \varphi_0^{-1}(\cup \mathcal{D}_0^i)$,
		\item \label{suru} if $i > 0$ then for each $C \in \mathfrak{C}^{m_{i-1}}$ there exists $C' \in \mathcal{D}^i_0$ with $C' \subseteq C$,
		\item \label{szorzas}
		
		for each $t\leq i$ and $j \in \omega$
		\[ \varphi_t\left[ \left( \bigcup_{k=0}^{i} \varphi_k^{-1}(\cup \mathcal{D}^i_j)\right) \cdot \left( \bigcup_{k=0}^{i} \varphi_k^{-1}( \cup \mathcal{D}^i_j) \right) \right]  \subseteq \cup \mathcal{D}^i_{j+1}, \]
		in other words, if $x \in \varphi_{r}^{-1}(\cup \mathcal{D}^i_j)$, $y \in \varphi_{s}^{-1}( \cup \mathcal{D}^i_j)$ ($r,s \leq i$), and $t \leq i$ is an integer for which $xy \in \dom(\varphi_t)$, then $\varphi_t(xy) \in \cup \mathcal{D}^i_{j+1}$,
		
		\item \label{inv} 
		
		for each $t \leq i$, $j \in \omega$ the following inclusion relation holds
		\[ \varphi_t\left[ \left( \bigcup_{k=0}^{i} \varphi_k^{-1}(\cup \mathcal{D}^i_j)\right)^{-1} \right]  \subseteq \cup \mathcal{D}^i_{j+1} \]
		in other words, if $r \leq i$, where $x \in \varphi_{r}^{-1}(\cup \mathcal{D}^i_j)$,  and $t \leq i$ is an integer for which $x^{-1} \in \dom(\varphi_t)$, then $\varphi_t(x^{-1}) \in \cup \mathcal{D}^i_{j+1}$,
		
		
		\item \label{mmertek} if $i>0$ then for every $r \leq i$ 
		\begin{equation} \label{mertek}
				\mu\left( \varphi_r^{-1}\left( \bigcup \mathcal{D}^i_i \right) \right) \leq \frac{1}{i^2} .
		\end{equation}
	\end{enumerate} 

	\bl \label{Dlemma}
			There is a strictly increasing sequence of non-negative integers $(m_i)_{i \in \omega}$ and a sequence $(\mathcal{D}^i_j)_{i \in \omega, j \in \omega}$  satisfying $\eqref{Dmon}- \eqref{mmertek}$,
	\el
	\bp
	    	Let 
	    	\begin{equation} \label{m0} 
	    	m_0 = 0, \end{equation} and 
	    	$\mathcal{D}^0_j = \{ [0,1]^d \} $ for every $j \in \omega$.
	    	Fix $i >0$.
	    	Assume that $m_1<m_2 < \dots < m_{i-1}$, and $\mathcal{D}^k_j \subseteq \mathfrak{C}^{m_k}$ ($k <i$, $j \in \omega$) are defined.
	    	
	    	Using $\eqref{egysegegyben}$ let $x_0 =\varphi_0(e) \in  [0,1]^d$, and define the finite set $A_0 \subseteq [0,1]^d$ by translating $x_0$ to each $C \in \mathfrak{C}^{m_{i-1}}$
	    	\begin{equation} \label{A0def}  A_0 = \left\{ x_0 + \left( \frac{s_1}{2^{m_{i-1}}}, \frac{s_2}{2^{m_{i-1}}}, \dots,  \frac{s_d}{2^{m_{i-1}}}\right) :  \ s_1,s_2, \dots, s_d \in \mathbb{Z} \right\} \cap [0,1]^d .  \end{equation}
	    	
	    	Define $B_0 = \bigcup_{j=0}^{i} \varphi_j^{-1}(A_0) \subseteq G$. Obviously 
	    	\begin{equation} \label{eAban} A_0 \ni x_0 = \varphi_0(e) \end{equation} 
	    	 implies that
	    	\begin{equation} \label{eB0}
		    	e \in B_0.
	    	\end{equation}
	    	
	    	Define the following operation $\mathcal{F}$ on $\mathcal{P}(G)$
	    	\[ \mathcal{F}: \mathcal{P}(G) \to \mathcal{P}(G), \]
	    	\begin{equation} \label{fop} \mathcal{F}(H) =(H \cup H^{-1}) \cdot (H \cup H^{-1}). \end{equation}
	    	Define inductively the  sequences $B_j,A_j$ ($j \leq i$) by the following procedure.
	    	If $B_{j-1}$, $A_{j-1}$ are given,
	    	then let
	    	\begin{equation} \label{Akov} A_j =  \bigcup_{k=0}^{i} \varphi_k[\mathcal{F}(B_{j-1})],  \end{equation}
	    	i.e. we choose $A_j \subseteq [0,1]^d$ by taking for each element of $\mathcal{F}(B_{j-1})$ all corresponding points  in the unit cube according to the $\varphi_k$-s ($k \leq i$). (Notice that there may be points in $\mathcal{F}(B_{j-1})$ which are not contained in $\dom(\varphi_k)$-s ($k \leq i$), thus there are no corresponding points in $A_j$.) 
	    	And $B_j$ will be the set of corresponding group elements according to the $\varphi_k$-s ($0 \leq k \leq i$), that is,
	    	\begin{equation} \label{Bdef} B_j = \bigcup_{k=0}^{i} \varphi_k^{-1}(A_j) .\end{equation}
	    	    	By induction on $j$, we prove that $e \in B_j$. By $\eqref{eB0}$,  $e \in B_0$, and if $e \in B_j$,
    	then by $\eqref{fop}$ $e \in \mathcal{F}(B_j)$, thus  $\eqref{Akov}$ and $e \in \dom(\varphi_0)$ by $\eqref{egysegegyben}$
	    yield that 	 $\varphi_0(e) \in A_{j+1}$.
	     This implies that 
	    	 \begin{equation} \label{etart} e \in \bigcup_{k=0}^{i}  \varphi_k^{-1}(A_{j+1}) = B_{j+1}. \end{equation} 
	    	 Now by $\eqref{etart}$ and $\eqref{fop}$ we get that 
	    	 \begin{equation} \label{BFB} B_j \subseteq \mathcal{F}(B_j). \end{equation} 
	    	 We now check that
	    	the $B_j$-s and the $A_j$-s form an increasing sequence. 
	    	Indeed, for each  $j \in \omega$ using $\eqref{Akov}$, $\eqref{BFB}$ and $\eqref{Bdef}$
	    	\[ A_{j+1} =  \bigcup_{k=0}^{i}\varphi_k[ \mathcal{F}(B_j) ]  \supseteq \bigcup_{k=0}^{i} \varphi_k[ B_j ] = \bigcup_{k=0}^{i} \left(  \varphi_k\left[ \bigcup_{l=0}^{i} \varphi_l^{-1}(A_j) \right] \right) \supseteq A_j \ , \]
	    	and from that
	    	\[ B_{j+1} = \bigcup_{k=0}^{i} \varphi_k^{-1}(A_{j+1}) \supseteq \bigcup_{k=0}^{i} \varphi_k^{-1}(A_{j}) = B_j \ .\]


	    	For each $k$ define
	    	\[ \mathcal{C}^k = \{ C \in \mathfrak{C}^k: \ A_i \cap C \neq \emptyset \}, \]
	    	i.e. $\mathcal{C}^k$ is a minimal, at most $2^d|A_i|$-many element subset of $\mathfrak{C}^k$ 
	    	such that 
	    	\begin{equation} \label{Attart} A_i  \subseteq \inte_{[0,1]^d}(\cup \mathcal{C}^k) \end{equation}

	         (where $\inte_{[0,1]^d}$ denotes the interior with respect to the subspace topology of the cube). 
	    	Clearly $\bigcup \mathcal{C}^k \supseteq \bigcup \mathcal{C}^{k+1}$, and  $\bigcap_{k=0}^{\infty}\bigcup \mathcal{C}^k = A_i$.
	    	This implies that, since $\varphi_r$-s are homeomorphisms 
	    	\[ \bigcap_{k=0}^{\infty} \varphi_r^{-1}(\bigcup \mathcal{C}^k) = \varphi_r^{-1}(A_i). \]
			This intersection is finite, hence is of measure zero since $G$ is of positive dimension. 
	    	Since compact sets have finite measure, and $\dom(\varphi_r)$ is compact, there exists a $k_r$, such that
	    	\begin{equation} \label{inegyzet} \mu\left( \varphi_r^{-1}\left( \bigcup \mathcal{C}^{k_r} \right) \right) \leq \frac{1}{i^2}. \end{equation}
	    	Define 
	    	\[ l_i = \max \{ m_{i-1}+1, k_0, k_1, \dots, k_i \}. \]
	    	Let $\mathcal{E}_i = \mathcal{C}^{l_i} \subseteq \mathfrak{C}^{l_i}$,
	    	thus  for each $r \leq i$
	    	\begin{equation} \label{Emertek} \mu\left( \varphi_r^{-1}\left( \bigcup \mathcal{E}_i \right) \right) \leq \frac{1}{i^2}. \end{equation}
	        Now beginning from $j= i$, and stepping down with $j$ to $0$, we will define a sequence $\mathcal{E}_j \subseteq \mathfrak{C}^{l_j}$ for some $l_j$ depending only on $j$ and the already defined $\mathcal{E}_k$-s ($ k > j$, $k \leq i$), for which the following holds
	    	\begin{itemize}
	    		\item the $l_j$-s are nonincreasing, i.e. \begin{equation} \label{lcsokk} l_0 \geq l_{1} \geq \dots \geq l_i, \end{equation}
	    		\item for each $j \leq i$\begin{equation} \label{ooo} \mathcal{E}_j \subseteq \mathfrak{C}^{l_j}, \end{equation}
	    		\item for each $j<i$ \begin{equation} \label{Emon} \cup \mathcal{E}_{j+1} \supseteq \cup \mathcal{E}_j, \end{equation}
	    		\item $A_j$ is in the relative interior (i.e. in the cube) of $\cup \mathcal{E}_j$:
	    		\begin{equation}  \label{lab}
	    		\inte_{[0,1]^d}(\cup \mathcal{E}_j) \supseteq A_j,
	    		\end{equation} 
	    		
	    		\item if (for $r,s \leq i$) $x \in \varphi_r^{-1}(\cup \mathcal{E}_j)$, $y \in \varphi_s^{-1}(\cup \mathcal{E}_j)$ and $t\leq i$ is such that $xy \in \dom(\varphi_k)$, then $\varphi_t(xy) \in \cup \mathcal{E}_{j+1}$, in other words
	    		\begin{equation} \label{masodik}
	    		(\forall r,s \leq i, \ \forall t \leq i) \ \ \varphi_t[ \varphi_r^{-1}(\cup \mathcal{E}_j) \cdot \varphi_s^{-1}( \cup \mathcal{E}_j)) ] 
	    		\subseteq \cup \mathcal{E}_{j+1}.
	    			    		\end{equation}
%
			   \item If $x \in \varphi_r^{-1}(\cup \mathcal{E}_j)$ ($r \leq i$), and $t \leq i$ is such that $x^{-1} \in \dom(\varphi_t)$
			   then $\varphi_t(x^{-1}) \in \cup \mathcal{E}_{j+1}$, i.e.:
			   	\begin{equation} \label{harm} (\forall r \leq i, \ \forall t \leq i) \ \ \varphi_t[ (\varphi_r^{-1}(\cup \mathcal{E}_j))^{-1}]  \subseteq \cup \mathcal{E}_{j+1}   \end{equation}
	   
	   	\end{itemize}
	   	Obviously (by $\eqref{Attart}$ and the definition of $\mathcal{E}_i$) $\eqref{lab}$ holds for $\mathcal{E}_i$.
	   	
	   	\bl \label{Elemma}
	   		If a given $\mathcal{E}_i \subseteq \mathfrak{C}^{l_i}$ satisfies $\eqref{lab}$, then there are sequences $l_0,l_1, \dots l_{i-1}$ and $\mathcal{E}_0, \mathcal{E}_1, \dots, \mathcal{E}_{i-1}$
	   		such that $\eqref{lcsokk}-\eqref{harm}$ hold.
	   	\el
	   	\bp
	   			 Suppose that $\mathcal{E}_{j+1} \subseteq \mathfrak{C}^{l_{j+1}}, \mathcal{E}_{j+2}\subseteq \mathfrak{C}^{l_{j+2}}, \dots, \mathcal{E}_{i}\subseteq \mathfrak{C}^{l_{i}}$ are already defined, and $A_{j+1} \subseteq \cup \inte_{[0,1]^d}(\mathcal{E}_{j+1})$, i.e. $\eqref{lab}$ holds for $\mathcal{E}_{j+1}$.   	
	   			 

	   			 For the construction of $\mathcal{E}_j$ we will need the following claims.
	   			 \bc\label{clegyes}
	   			 	 For each $x \in A_j$ there are sufficiently small neighborhoods $U_x$, $V_x$ in the cube so that whenever
	   			 	 $y,z \in A_j$, $r,s \leq i$, then  for any $t \leq i$
	   			 	 \begin{equation} \label{UVsz} \varphi_t[ \varphi_r^{-1}(U_z) \cdot \varphi_s^{-1}(V_y)] \subseteq \cup \mathcal{E}_{j+1}, \end{equation}
	   			 	 i.e.
	   			 	 \[ p \in \varphi_r^{-1}(U_z) \cdot \varphi_s^{-1}(V_y) \cap \dom(\varphi_t) \ \Rightarrow \ \varphi_t(p) \in \cup \mathcal{E}_{j+1} .\]
	   			 \ec 
	   			 \bc\label{clW}
	   			 	 For any $x \in A_j$ there is a neighborhood $W_x \subseteq [0,1]^d$ such that 
	   			 	 for each $r,t \leq i$
	   			 	 \begin{equation} \label{WW}
	   			 	 \varphi_t[(\varphi_r^{-1}(W_x))^{-1}] \subseteq \cup \mathcal{E}_{j+1}. 
	   			 	 \end{equation}
	   			 \ec
	   			 \bp(Claim $\ref{clegyes}$)
	   			 	 For defining the $U_x$-s and $V_x$-s we will need for any fixed $r,s,t \leq i$, and $y,z \in A_j$  neighborhoods $U^t_{r,y,s,z} $ of $y$, and $ V^t_{r,y,s,z}$ of $z$ 
	   			 	 such that
	   			 	 \begin{equation} \label{clai}
	   			 	 \varphi_t[ \varphi_r^{-1}(U^t_{r,y,s,z}) \varphi_s^{-1}(V^t_{r,y,s,z}) ] \subseteq \cup \mathcal{E}_{j+1} 
	   			 	 \end{equation}
	   			 	 holds.
	   			 	 Let $y,z \in A_j$, $r,s,t \leq i$  be fixed. 	Then, if $g =\varphi_r^{-1}(y)$, $h =\varphi_s^{-1}(z)$,
	   			 	 then $g,h \in B_j$ (by the definition of the $B_k$-s, $\eqref{Bdef}$), $gh \in \mathcal{F}(B_j)$ (by  $\eqref{fop}$).
	   			 	 Now we have two cases depending on whether $gh \in \dom(\varphi_t)$ holds.

	   			 	 First if $t \leq i$ is such that $gh \in \dom(\varphi_t)$:
	   			 	 	then first, $a =\varphi_t(gh) \in A_{j+1}$ (by $\eqref{Akov}$). 
	   			 	 	Since $\mathcal{E}_{j+1}$ was choosen so that $A_{j+1} \subseteq \inte_{[0,1]^d}(\cup \mathcal{E}_{j+1})$ by $\eqref{lab}$,
	   			 	 	$a \in \inte_{[0,1]^d}(\cup \mathcal{E}_{j+1})$.
	   			 	 	Moreover $\varphi_t$ is a homeomorphism,  thus $\varphi_t^{-1}(\inte_{[0,1]^d}(\cup \mathcal{E}_{j+1}))$ is relatively open in $\dom(\varphi_t)$.
	   			 	 	Therefore there is an open set $W$ in $G$ containing $gh$, such that 
	   			 	 	\begin{equation} \label{Wdef}
	   			 	 	W \cap \dom(\varphi_t) \subseteq \varphi_t^{-1}(\inte_{[0,1]^d}(\cup \mathcal{E}_{j+1})) . 
	   			 	 	\end{equation}
	   			 	  	Then, by the continuity of the multiplication, we have that there are
	   			 	 	open sets  $U,V \subseteq G$, $g \in U$,  $h \in V$,
	   			 	 	for which 
	   			 	 	\begin{equation}  \label{UVdef} UV \subseteq W.\end{equation}
	   			 	 	Applying $\varphi_t$ to both sides, and using $\eqref{Wdef}$ we get that
	   			 	 	\begin{equation} \label{Gben} 
	   			 	 	\varphi_t \left[  U \cdot V  \right] \subseteq \varphi_t[\varphi_t^{-1}(\inte_{[0,1]^d}(\cup \mathcal{E}_{j+1}))] =  \inte_{[0,1]^d}(\cup \mathcal{E}_{j+1}). \end{equation}
	   			 	 	Let $U^t_{r,y,s,z} = \varphi_r[U] = \varphi_r[U \cap \dom(\varphi_r) ] \ni y$, $V^t_{r,y,s,z} = \varphi_s[V] = \varphi_s[V \cap \dom(\varphi_s)] \ni z$  relatively open sets in the cube ($\varphi_s$, $\varphi_r$ are homeomorphisms), then $\eqref{Gben}$ implies that
	   			 	 	\[ \varphi_t\left[ \varphi_r^{-1}(U^t_{r,y,s,z}) \varphi_s^{-1}(V^t_{r,y,s,z}) \right] \subseteq  \inte_{[0,1]^d}(\cup \mathcal{E}_{j+1}). \]

	   			 	Second, if $t$ is such that $gh \notin \dom(\varphi_t)$,
	   			 	 	then (recall that $\varphi_t$ is a homeomorphism, thus $\varphi_t^{-1}([0,1]^d)$ is compact),
	   			 	 	there are open sets $U \ni g$, $V \ni h$ in $G$
	   			 	 	such that 
	   			 	 	\begin{equation} \label{diszjj} UV \cap \dom(\varphi_t) = \emptyset .\end{equation}
	   			 	 	Let  $U^t_{r,y,s,z} = \varphi_r[U] \ni y$, $V^t_{r,y,s,z} = \varphi_s[V] \ni z$ be open sets in the cube, then $\eqref{diszjj}$ clearly yields that
	   			 	 	\[ \varphi_r^{-1}(U^t_{r,y,s,z}) \varphi_s^{-1}(V^t_{r,y,s,z}) \cap \dom(\varphi_t) = \emptyset, \]
	   			 	 	thus
	   			 	 	\[  \varphi_t\left[ \varphi_r^{-1}(U^t_{r,y,s,z}) \varphi_s^{-1}(V^t_{r,y,s,z}) \right] = \emptyset \subseteq  \inte_{[0,1]^d}(\cup \mathcal{E}_{j+1}) . \]

	   			 	 We obtain that in both cases  	$\eqref{clai}$ holds for the  neighborhoods $U^t_{r,y,s,z} $ of $y$, and $ V^t_{r,y,s,z}$ of $z$.
	   			 	 
	   			 	 Now we can define the $U_x$-s and $V_x$-s ($x \in A_j$) as
	   			 	 \begin{equation} \label{Udef} 
	   			 	 U_x = \bigcap_{s,r, t \leq i, z \in A_j} U^t_{r,x,s,z} ,
	   			 	 \end{equation}
	   			 	 and
	   			 	 \begin{equation} \label{Vdef}
	   			 	 V_x = \bigcap_{s,r,t \leq i, y \in A_j} V^t_{r,y,s,x}.
	   			 	 \end{equation}	 
	   			 	 Then $U_x, V_x \ni x$ are relatively open sets in the cube (since $A_j$ is finite), and clearly for  arbitrary 
	   			 	 $r,s,t \leq i$, and $y,z \in A_j$ by $\eqref{Udef}$, $\eqref{Vdef}$ and $\eqref{clai}$
	   			 	 \[ \varphi_t[\varphi_r^{-1}(U_y) \cdot \varphi_{s}^{-1}(V_z)] \subseteq  \varphi_t[ \varphi_r^{-1}(U^t_{r,y,s,z}) \cdot   \varphi^{-1}_{s}(V^t_{r,y,s,z}) ]  \subseteq  \cup \mathcal{E}_{j+1},
	   			 	 \]
	   			 	 
	   			 	 i.e. $\eqref{UVsz}$ holds, indeed.
	   			 \ep

	   			 
	   			 We can turn to the proof of our second claim.
	   			\bp (Claim $\ref{clW}$)

	   			 If $x \in A_j$ is fixed, for the desired neighborhood $W_x$ we do the following.	   			
	   			 First	for each $x \in A_j$ and $r,t \leq i$ we would like to define a neighborhood $W^t_{x,r}$ of $x$ in $[0,1]^d$ such that 
	   			 \begin{equation} \label{Wclaim} \varphi_t[(\varphi_r^{-1}(W^t_{x,r}))^{-1}] \subseteq \cup \mathcal{E}_{j+1} . \end{equation} 
	   			 Let $x \in A_j$ and $r,t \leq i$ be fixed. Set $g= \varphi_r^{-1}(x)$.
	   			
	   			If $g^{-1} = (\varphi_r^{-1}(x))^{-1} \in \dom(\varphi_t)$,    			 	then
	   			 	since $g \in B_j$ by $\eqref{Bdef}$, thus $g^{-1} \in \mathcal{F}(B_{j})$ (by $\eqref{fop}$, and $e \in B_j$).  
	   			 	By $\eqref{Akov}$ $\varphi_t(g^{-1}) \in A_{j+1}$.
	   			 	We know that $A_{j+1} \subseteq \inte_{[0,1]^d}(\cup \mathcal{E}_{j+1})$, i.e. $\varphi_t(g^{-1})$ is an element of the open set   
	   			 	$\inte_{[0,1]^d}(\cup \mathcal{E}_{j+1}) \subseteq [0,1]^d$,  hence $g^{-1}$ is the element of the relatively open set 
	   			 	$\varphi_t^{-1}(\inte_{[0,1]^d}(\cup \mathcal{E}_{j+1})) \subseteq \dom(\varphi_t)$.
	   			 	From which there is an open $V \subseteq G$ containing $g^{-1}$ such that 
	   			 	\[ V \cap \dom(\varphi_t) \subseteq \varphi_t^{-1}(\inte_{[0,1]^d}(\cup \mathcal{E}_{j+1})). \]
	   			 	Now $V^{-1} \subseteq G$ is an open neighborhood of $g$.
	   			 	If we set $W_{x,r}^t = \varphi_r[V^{-1}] = \varphi_r[V^{-1} \cap \dom(\varphi_r)]$, then $x \in W_{x,r}^t$ and
	   			 	\[ (\varphi_r^{-1}(W_{x,r}^t) )^{-1} \subseteq (V^{-1})^{-1} = V, \]
	   			 	from which it follows that
	   			 	\[ (\varphi_r^{-1}(W_{x,r}^t) )^{-1} \cap \dom(\varphi_t) \subseteq V \cap \dom(\varphi_t) \subseteq  \varphi_t^{-1}(\inte_{[0,1]^d}(\cup \mathcal{E}_{j+1})), \]
	   			 	i.e. 
	   			 	\[ \varphi_t[(\varphi_r^{-1}(W_{x,r}^t) )^{-1}] \subseteq \inte_{[0,1]^d}(\cup \mathcal{E}_{j+1}) \subseteq \cup \mathcal{E}_{j+1} .\]
	   		
	   			 In the other case, when $g^{-1} = (\varphi_r^{-1}(x))^{-1} \notin \dom(\varphi_t)$, then by the compactness of $\dom(\varphi_t)$, there is an open 
	   			 	neighborhood $V$ of $g^{-1}$ such that $V \cap \dom(\varphi_t ) = \emptyset$. But then 
	   			 	 $V^{-1}$ is a neighborhood of $g$ such that $(V^{-1})^{-1} = V \subseteq G \setminus \dom(\varphi_t)$.
	   			 	Then, with $W_{x,r}^t = \varphi_r[V^{-1}] = \varphi_r[V^{-1} \cap \dom(\varphi_r)]$ clearly $x \in W_{x,r}^t$ and
	   			 	 $(\varphi_r^{-1}(W^t_{x,r}))^{-1} \subseteq V$ is disjoint from $\dom(\varphi_t)$,
	   			 	thus 
	   			 	\[ \varphi_t[(\varphi_r^{-1}(W_{x,r}^t))^{-1}] = \emptyset \subseteq \cup \mathcal{E}_{j+1}. \]

	   			 We conclude that in both cases i.e. independently of whether $g^{-1} \in \dom(\varphi_t)$ $\eqref{Wclaim}$ holds.
	   			 
	   			 Now we can define the $W_x$-s ($x \in A_j$)
	   			 \begin{equation} \label{WWdef} W_x =  \bigcap_{r, t \leq i} W^t_{x,r} \ .\end{equation}
	   			 Then, fixing $x \in A_j$, for any $r,t \leq i$, using $\eqref{WWdef}$ and $\eqref{Wclaim}$
	   			 \[ \varphi_t[(\varphi_r^{-1}(W_x))^{-1}] \subseteq \varphi_t[(\varphi_r^{-1}(W^t_{x,r}))^{-1}] \subseteq  \cup \mathcal{E}_{j+1} , \]
	   			 i.e. $\eqref{WW}$ holds, indeed.

	   			 \ep
	   			 
	   			 %

	   			 Now having these claims proven, we are ready to finish the proof of Lemma $\ref{Elemma}$.
	   			 We have that $A_j$ is contained in the open sets  $\bigcup_{x \in A_j} U_x$, $ \bigcup_{x \in A_j} V_x$ and $\bigcup_{x \in A_j} W_x$. 
	   			 For each $p \geq l_{j+1}$ let $\mathcal{C}^p_j \subseteq \mathfrak{C}^p$ denote
	   			 \begin{equation}  \mathcal{C}^p_j = \{ C \in \mathfrak{C}^p: \ A_j \cap C \neq \emptyset \}, \end{equation}
	   			  so this is the unique minimal covering of $A_j$, for which $A_j \subseteq \inte_{[0,1]^d}(\cup \mathcal{C}^p_j)$ (hence $|\mathcal{C}^p_j| \leq 2^d |A_j|$).
	   			   Furthermore, since  $A_j$-s form an increasing sequence, i.e. $A_j \subseteq A_{j+1}$ (by $\eqref{Amon}$), combining $\eqref{lab}$, the inequality $p \geq l_{j+1}$, and $\mathcal{E}_{j+1} \subseteq \mathfrak{C}^{l_{j+1}}$ (by $\eqref{ooo}$)
	   			   \begin{equation} \label{nov} \cup \mathcal{C}^p_j \subseteq \cup \mathcal{E}_{j+1}. \end{equation}
	   			  
	   			 Now if $p \geq l_{j+1}$ is sufficiently large, then
	   			 \begin{equation} \label{fet} A_j \subseteq \cup \mathcal{C}^p_j \subseteq  \left( \bigcup_{x \in A_j} U_x \right) \cap \left(\bigcup_{y \in A_j} V_y \right) \cap \left( \bigcup_{x \in A_j} W_x \right). \end{equation}
	   			  
	   			 Fix such a $p \geq l_{j+1}$ for which
	   			 $\eqref{fet}$ and $\eqref{nov}$ hold, and let  $\mathcal{E}_j = \mathcal{C}^p_j$, $l_j = p$. 
	   			 ($\eqref{lcsokk}$ obviously holds.)
	   			 
	   			 This clearly implies that $\eqref{Emon}$ and $\eqref{lab}$ hold,   we only have to check that $\eqref{masodik}$ and $\eqref{harm}$ hold.
	   			 First we verify that $\eqref{masodik}$ holds, i.e. for any $t \leq i$, $r,s \leq i$
	   			 \[ \varphi_t\left[ \varphi_r^{-1}( \cup \mathcal{E}_j)  \varphi_s^{-1}( \cup \mathcal{E}_j)  \right] \subseteq  \cup \mathcal{E}_{j+1} .\]
	   			 Fixing such $t$, $r$, $s$, since  $\cup \mathcal{E}_j \subseteq \left( \bigcup_{x \in A_j} U_x \right)$, and $\cup \mathcal{E}_j \subseteq \left(\bigcup_{y \in A_j} V_y \right)$ (by $\eqref{fet}$ and the definition of $\mathcal{E}_j$)
	   			 \[	\varphi_t\left[  \varphi_r^{-1}( \cup \mathcal{E}_j) \varphi_s^{-1}( \cup \mathcal{E}_j)  \right] \subseteq \varphi_t\left[  \varphi_r^{-1} \left( \bigcup_{x \in A_j} U_x \right) \varphi_s^{-1} \left(\bigcup_{y \in A_j} V_y \right)  \right], \]
	   			 and using $\eqref{UVsz}$
	   			 \[ \varphi_t\left[  \varphi_r^{-1} \left( \bigcup_{x \in A_j} U_x \right) \varphi_s^{-1} \left(\bigcup_{y \in A_j} V_y \right)  \right] \subseteq \bigcup_{x,y \in A_j} \varphi_t\left[  \varphi_r^{-1} ( U_x) \varphi_s^{-1}  (V_y )  \right] \subseteq \cup \mathcal{E}_{j+1},\]
	   			 i.e. $\eqref{masodik}$ holds, indeed.
	   			 
	   			 Similarly,  since  $\cup \mathcal{E}_j \subseteq \left( \bigcup_{x \in A_j} W_x \right)$, for any fixed $r,t \leq i$

	   			 by $\eqref{fet}$ (and $\mathcal{C}^{l_j} = \mathcal{E}_j$)
	   			 \[ \left( \varphi_r^{-1}( \cup \mathcal{E}_j) \right)^{-1} \subseteq  \left( \varphi_r^{-1} \left( \bigcup_{x \in A_j} W_x \right) \right)^{-1}  
	   			 	   			 =   \bigcup_{x \in A_j} (\varphi_r^{-1}(W_x))^{-1}, \]
	   			 and  applying $\varphi_t$ for this equation, $\eqref{WW}$ gives that
	   			 \[  \varphi_t\left[ (  \varphi_r^{-1}( \cup \mathcal{E}_j) )^{-1} \right] \subseteq  \varphi_t \left[ \bigcup_{x \in A_j} (\varphi_r^{-1}(W_x))^{-1} \right] \subseteq (\cup \mathcal{E}_{j+1}) . \]
	   			 This ends the proof of Lemma $\ref{Elemma}$.
     \ep

	    	Now we finish the proof of Lemma $\ref{Dlemma}$.
	    	
	    	At the point when $\mathcal{E}_j$ is already defined for all $0 \leq j \leq i$, and $l_0$ denotes the fineness of $\mathcal{E}_0$, i.e. 
	    	$\mathcal{E}_0 \subseteq \mathfrak{C}^{l_0}$, then we define $m_i$ and the $\mathcal{D}_j^i$-s ($j \in \omega$) as follows.
	    	Let $m_i = l_0$, and for each $j \leq i$ 
	    	let $\mathcal{D}_j^i \subseteq \mathfrak{C}^{m_i}$ be such that 
	    	\begin{equation} \label{DdfE} \cup \mathcal{D}^i_j = \cup \mathcal{E}_j, \ \ (j \leq i) \end{equation}
	    	 i.e. a subdivision of $\mathcal{E}_j$ (which is possible since $\mathcal{E}_j \subseteq \mathfrak{C}^{l_j}$, where $l_j$-s form a decreasing sequence by $\eqref{lcsokk}$, and $m_i = l_0$).
	    		For $j > i$, let 
	    		\begin{equation} \label{Ddf2} \mathcal{D}^i_j = \mathfrak{C}^{m_i}, \ \ (j > i) \end{equation}
	    		 i.e. we choose all the $2^{m_i}$-many small cubes.
	    		We only have to check that $\eqref{Dmon}- \eqref{mmertek}$ hold.
	    	
	    	$\eqref{Emon}$ implies that $\eqref{Dmon}$ holds for $j < i$, and clearly $\eqref{Ddf2}$ implies for $j \geq i$.
	    	Combining $\eqref{eAban}$ with $\eqref{lab}$ one can get that $e \in \cup \mathcal{E}_0$ thus by $\eqref{DdfE}$ we obtain $\eqref{De}$.
	    	
	    		Now $\inte_{[0,1]^d}(\cup \mathcal{E}_0) \supseteq A_0$ (by $\eqref{lab}$) and $A_0$ contains at least one point from each $C \in \mathfrak{C}^{m_{i-1}}$ (by $\eqref{A0def}$),
	    		therefore for each $C \in \mathfrak{C}^{m_{i-1}}$ $\mathcal{E}_0$ contains a cube $C'$ with $C' \subseteq C$. 
	    		    		    	Hence from the definition of $\mathcal{D}^i_0$ (i.e. $\cup \mathcal{D}^i_0 = \cup \mathcal{E}_0$)
	    		    		    	$\mathcal{D}^i_0$ also contains a cube $C''$ with $C'' \subseteq C$. This means that $\eqref{suru}$ holds.
	    	
	    	Since $\eqref{masodik}$ and $\eqref{harm}$ hold, and because $\cup \mathcal{E}_j = \cup \mathcal{D}^i_j$ ($j \leq i$), $\eqref{szorzas}$ and $\eqref{inv}$ obviously hold for the $\mathcal{D}^i_j$-s where $j<i$. If $j \geq i$ then $\eqref{Ddf2}$ obviously implies that $\cup \mathcal{D}^i_{j+1} = [0,1]^d$, therefore we obtain $\eqref{szorzas}$ and $\eqref{inv}$.
	    	
	    	Finally $\cup \mathcal{E}_i = \cup \mathcal{D}_i^i$ and $\eqref{Emertek}$ yield equation $\eqref{mertek}$.

	\ep

     Now we return to the proof of Proposition $\ref{Lieeset}$, and  construct the desired subgroup.
     
     For each $z \in [0,1]^d$ let $T_z$ denote the set of the those sequences of $\prod_{i \in \omega} \mathfrak{C}^{m_i}$ whose intersection is $\{z\}$, i.e.
     \[ T_z = \left\{ (C_i)_{i \in \omega}| \ \ \forall i C_i \in \mathfrak{C}^{m_i}, \ \cap_{i \in \omega} C_i = \{ z \}   \right\}.\]
     Now choose a non-principal ultrafilter $\mathcal{U}$ on $\omega$, and define $S' \subseteq [0,1]^d$ in the following way:
     \[ S' = \left\{ z : \ \  \exists l \ \exists (C_i)_{i \in \omega} \in T_z \  \{i\in \omega| \   C_i \in \mathcal{D}^i_l \} \in \mathcal{U}   \right\}  \]
     Let 
     \begin{equation} \label{Sdefje} S = \bigcup_{i=0}^{\infty} \varphi_i^{-1}(S'). \end{equation}
     Next we will show that $S$ is a subgroup of $G$, $S$ is null, and is not meager:
     \begin{enumerate}[(i)]
     	\item $S$ is a subgroup.
	\end{enumerate}    
      First, $S$ is closed under the multiplication.
    	
	Let $x,y \in S$	and $k_1,k_2,k_3 \in \omega$ be such that $x \in \varphi^{-1}_{k_1}(S')$, $y \in \varphi^{-1}_{k_2}(S')$,
     		$xy \in \dom(\varphi_{k_3})$, 
     		$\varphi_{k_1}(x) = x',\varphi_{k_2}(y) = y' \in S'$.
     		Let $(C_{i}^x)_{i  \in \omega} \in T_{x'}$, 
     		$(C_{i}^y)_{i \in \omega} \in T_{y'}$,
     		$U_x, U_y \in \mathcal{U}$, $l_x$,$l_y$ be such that
     		\begin{equation} \label{xes} \{i : \ C_{i}^x \in \mathcal{D}^i_{l_x} \} = U_x,  \end{equation}
     		\begin{equation} \label{yos} \{ i : \ C_{i}^y \in \mathcal{D}^i_{l_y} \} = U_y. \end{equation}
     		Recall that $\cup \mathcal{D}^i_k \subseteq \cup \mathcal{D}^i_{k+1}$ ($k \in \omega$) by the construction of the $\mathcal{D}^i_k$-s, see $\eqref{Dmon}$.
     			Let $l = \max \{l_x,l_y \}$ denote the larger one, then
     	 if  $i \in U_x \cap U_y \in \mathcal{U}$,  using $\eqref{xes}$, $\eqref{yos}$
     		 \[ x' = \varphi_{k_1}(x) \in C^x_i \in \mathcal{D}^i_{l_x}  \subseteq \mathcal{D}^i_l \]
     		 and
     		 \[ y' = \varphi_{k_2}(y) \in C^y_i \in \mathcal{D}^i_{l_y}  \subseteq \mathcal{D}^i_l, \]
     		 i.e.
     	
     			\begin{equation} \label{xes2} \{i : \ C^x_i \in \mathcal{D}^i_{l} \} \supseteq  U_x  \supseteq U_x \cap U_y \in \mathcal{U},  \end{equation}
     			\begin{equation} \label{yos2} \{ i : \ C^y_i \in \mathcal{D}^i_{l} \} \supseteq U_y \supseteq U_x \cap U_y \in \mathcal{U} . \end{equation}

     		From which,
     		  \[ i \in U_x \cap U_y \ \to \ xy \in \left( \varphi_{k_1}^{-1}(\cup \mathcal{D}^i_l) \cdot \varphi_{k_2}^{-1}(\cup \mathcal{D}^i_l) \right) \cap \dom(\varphi_{k_3}). \]
     		Now $i$ is such that $i \geq \max \{ k_1,k_2,k_3 \} =:k$, we can use the property $\eqref{szorzas}$ of the $\mathcal{D}^i_j$-s, that implies $xy \in \varphi_{k_3}^{-1}(\cup \mathcal{D}^i_{l+1})$, i.e. for each $i \in U_x \cap U_y \cap [k, \infty)$ there is an $C_i' \in \mathcal{D}^i_{l+1}$, for which
     		\[ \varphi_{k_3}(xy) \in C_i' .\]
     		This is true for all $i \in U_x \cap U_y \in \mathcal{U}$, $i \geq k$, thus there is a partial sequence $(C_i')_{i \in U_x \cap U_y \cap [k,\infty)}$ for which
     		\[  \varphi_{k_3}(xy) \in C_i'  \ \ (i \in  U_x \cap U_y \cap [k, \infty)).\] 
     		Choosing  $C_i'$-s for the remaining $i$-s ( $\omega \setminus (U_x \cap U_y \cap [k, \infty) )$ ) so that $\varphi_{k_3}(xy) \in C_i'$ hold will  yield an appropriate sequence 
     		$(C_i')_{i \in \omega}$ from $T_{\varphi_{k_3}(xy)}$, i.e.
     		\[ \{\varphi_{k_3}(xy) \} = \bigcap_{i \in \omega} C_i', \]
     		and
     		\[ \{ i \in \omega: \  C_i' \in \mathcal{D}^i_{l+1} \} \supseteq U_x \cap U_y \cap [l, \infty) \in \mathcal{U}. \]
     		Hence  $(C_i')_{i \in \omega}$, $l+1$, $U_x \cap U_y \cap [k, \infty)]$ witness that $\varphi_{k_3}(xy) \in S'$, 
     		 we can conclude that 
     		\[  xy \in S .\]

     	Similarly, for the closedness under taking inverse
     	  	fix $x \in S$, where $\varphi_{k_1}(x) = x' \in S'$ is witnessed by the aforementioned ($(C_i)_{i \in \omega}$, $l_x$, $U_x$).
     		 Now, if $i \in U_x$:
     		 \[ \varphi_{k_1}(x) = x' \in C_i \in \mathcal{D}^i_{l_x}. \]
     		 Fix $k$  such that $x^{-1} \in \dom(\varphi_k)$. Now if $i \in U_x \cap [k_x, \infty) \cap [k, \infty)$ then $\eqref{inv}$ implies
     		 \[ \varphi_{k}(x^{-1}) \in \cup \mathcal{D}^i_{l_x+1}, \]
     		 so there is an $C_i' \in \mathcal{D}^i_{l_x+1}$ for which $\varphi_k(x^{-1}) \in C_i'$. Again, this is true for each $i \in U_x \cap [\max\{k_x,k \}, \infty)$, i.e. we have a sequence $(C_i')_{i \in \omega}$ from $T_{\varphi_k(x^{-1})}$, such that for all $i \in U_x$ if $i \geq k,k_x$, then $C_i' \in \mathcal{D}^i_{l_x+1}$. Since $U_x \cap [k, \infty) \cap [k_x, \infty) \in \mathcal{U}$, this verifies that $\varphi_k(x^{-1}) \in S'$, hence $x^{-1} \in S$.

     	\begin{enumerate}[(ii)]	
     		\item $S$ is null.
         \end{enumerate} 
     	
     	By the $\sigma$-additivity, it is enough to show that $\mu( \varphi_p^{-1}(S') )= 0$ for each $p \in \omega$.
     	  And since 
     	\[ S' = \bigcup_{l \in \omega} \{ z \in [0,1]^d| \ \exists (C_{i})_{i \in \omega} \in T_z  \ \{i \in \omega| \ C_{i} \in \mathcal{D}_l^i \} \in \mathcal{U} \} \]
     	 it is sufficient to show that for each $l,p \in \omega$
     	\[ \mu \left(  \varphi_p^{-1} \left(  \{ z \in [0,1]^d | \ \exists (C_{i})_{i \in \omega} \in T_z  \ \{i \in \omega |  \ C_{i} \in \mathcal{D}_l^i \} \in \mathcal{U} \} \right) \right) = 0 .\]
     	But  
     	\[  \begin{array}{rrrl} \left\{ z \in [0,1]^d| \ \right. & \ \exists (C_{i})_{i \in \omega} \in T_z  \  & \{ i \in \omega | \ C_i \in \mathcal{D}_l^i \} \in \mathcal{U} \left. \right\} & \subseteq   \\
     		\left\{ z \in [0,1]^d | \ \right.  & \ \exists (C_{i})_{i \in \omega} \in T_z  \ & \{ i \in \omega | \ C_{i} \in \mathcal{D}_i^i \} \in \mathcal{U} \left. \right\} & \subseteq  \end{array} \]
     	since for fixed $l$, $\mathcal{D}_l^i \subseteq \mathcal{D}_i^i$ is true for all but finitely many $i$-s, and $\mathcal{U}$ is non-principal. Moreover, using the nonprincipality of $\mathcal{U}$ again (i.e. $\mathcal{U}$ only contains infinite subsets of $\omega$) the latter set is contained in a larger one
     	\[ Z:= \left\{ z \in [0,1]^d | \ \right.   \ \exists (C_{i})_{i \in \omega} \in T_z  \  \{ i \in \omega | \ C_{i} \in \mathcal{D}_i^i \} \in \mathcal{U} \left. \right\} \subseteq \] 
	     \[	     \left\{  z \in [0,1]^d |  \ \right. \ \exists (C_{i})_{i \in \omega} \in T_z  \  |\{ i \in \omega | \ C_{i} \in \mathcal{D}_i^i \}| = \infty \left. \right\}. 
   \]
   Then, if $z$ and $(C_{i})_{i \in \omega}$ are such that $\cap_{i \in \omega}C_i = \{z\}$ (i.e. $(C_{i})_{i \in \omega} \in T_z$), 
   and  $C_{i} \in \mathcal{D}_i^i $ for infinitely many $i$-s, what obviously implies that (for this fixed $z$) for infinitely many $i$-s there are sets $C_i' \ni z$,  $C'_{i} \in \mathcal{D}_i^i$. Therefore
    \[ 	Z \subseteq \bigcap_{k=1}^{\infty} \bigcup_{n=k}^{\infty} \{ z \in [0,1]^d |  \ \exists C'_{n} \ni z:   \ C'_n \in \mathcal{D}_n^n \}  = \]
     	\[
     	\bigcap_{k=1}^{\infty} \bigcup_{n=k}^{\infty} \ \cup \mathcal{D}_n^n . 
        	\]
     	Considering its preimage under the homeomorphism $\varphi_p$ we obtain that for any $r \in \omega$
     	\[ \varphi_p^{-1} \left( \bigcap_{k=1}^{\infty} \bigcup_{n=k}^{\infty} \ \cup \mathcal{D}_n^n \right) \subseteq \varphi_p^{-1}\left( \bigcup_{n=r}^{\infty} \ \cup \mathcal{D}_n^n\right), \]
     	and applying $\eqref{mertek}$  (i.e. $t \geq p$ implies that $\mu(\varphi_p^{-1} (\cup \mathcal{D}_t^t)) \leq \frac{1}{t^2}$), for a fixed $r \geq p$,
     	\[ \mu \left( \varphi_p^{-1}( \left\{ z \in [0,1]^d | \ \exists (C_i)_{i \in \omega} \in T_z  : \{i \in \omega | \ C_i \in \mathcal{D}_l^i \} \in \mathcal{U} \  \right\} ) \right) \leq \]
     	\[ \leq \mu\left( \varphi_p^{-1}\left(\bigcap_{k=1}^{\infty} \bigcup_{n=k}^{\infty} \ \cup \mathcal{D}_n^n \right)\right) \leq \mu\left( \varphi_p^{-1}\left( \bigcup_{n=r}^{\infty} \ \cup \mathcal{D}_n^n\right) \right) \leq \sum_{n=r}^\infty \frac{1}{n^2} \]
     	which upper bound tends to $0$ as $r$ tends to $\infty$.
     	
     	\begin{enumerate}[(iii)]
     	\item $S$ is not meager.
     	\end{enumerate}
     	     	
     	Let $R \subseteq G$ be co-meager, then it can be easily seen that $R \cap \dom(\varphi_0)$ is co-meager in $\dom(\varphi_0)$.
By $\eqref{Sdefje}$ it suffices to show that $\varphi_0[R] \cap S' \neq \emptyset$.
     	
     	Since $\varphi_0$ is a homeomorphism, the image  $\varphi_0[R] \subseteq [0,1]^d$ of $R \cap \dom(\varphi_0)$ is co-meager.
       We will apply  Lemma $\ref{rezid}$ in $[0,1]^d$ , for the co-meager set $\varphi_0[R]$, with the compact sets in $\bigcup_{i \in \omega} \mathfrak{C}^{m_i}$ in the following way. First, recall that for $k \in \omega$ $\mathfrak{C}^k$ consists of the $2^{dk}$-many cubes of sidelength $\frac{1}{2^k}$ with (pairwise disjoint interior and)  $\cup \mathfrak{C}^k = [0,1]^d$. In particular 
       $\mathfrak{C}^{m_0} = \mathfrak{C}^0 = \{ [0,1]^d \}$ by $\eqref{m0}$. Define $M_0 = \{ 1\}$, and
       $M_{i+1} = \{1,2, \dots, 2^{(m_{i+1}-m_i)d}  \}$. Let $C_1 = [0,1]^d \in \mathfrak{C}^{m_0}$, and suppose that the compact sets 
       $(C_{p_0p_1 \dots p_j})_{p_0p_1 \dots p_j \in \prod_{k \leq j} M_k}$ ($j \leq i$) are defined and satisfy the conditions in Lemma $\ref{rezid}$ (where $C_{p_0p_1 \dots p_j} \in \mathfrak{C}^{m_j}$). Then for a fixed set
       $C_{p_0p_1 \dots p_i} \in \mathfrak{C}^{m_i}$ as 
       \[ |M_{i+1}| = 2^ {(m_{i+1}-m_i)d} = |\{ C \in \mathfrak{C}^{m_{i+1}} \ | \ C \subseteq C_{p_0p_1 \dots p_i}  \}| \]
       	we can  define the $C_{p_0p_1 \dots p_ik}$-s ($k \in M_{i+1}$) to be different elements of the set $\{ C \in \mathfrak{C}^{m_{i+1}} \ | \ C \subseteq C_{p_0p_1 \dots p_i}  \}$. This implies that the assumptions of Lemma $\ref{rezid}$ hold.
	 Note that for each $i \in \omega$
	 \begin{equation} \label{frakC}
	 \{ C_{p_0p_1\dots p_i} \ | \ \ p_0p_1\dots p_i \in \prod_{j \leq i} M_j \} = \mathfrak{C}^{m_i},
	 \end{equation}
	 and for each $p_0p_1\dots p_i \in \prod_{j\leq i} M_j$
	 \begin{equation} \label{kitolt}
	  \bigcup_{k \in M_{i+1}} C_{p_0p_1 \dots p_i k} = C_{p_0p_1 \dots p_i}.
	 \end{equation}

       	The lemma gives a strictly increasing sequence $(n_j)_{j \in \omega}$, and a sequence $(r_j)_{j \in \omega}$,
     	 ($r_j \in M_j$), such that for any infinite sequence $(s_j)_{j \in \omega}$ ($s_j \in M_j$)
     	 \begin{equation} \label{kockaslemma}
     	    \text{if the set } \{l: r_{|[n_l,n_{l+1})} =s_{|n_l, n_{l+1}) }    \} \text{ is infinite} \ \Rightarrow \ \bigcap_{j \in \omega}  C_{s_0s_1 \dots s_j} \subseteq \varphi_0[R] . 
     	 \end{equation}
     	
     	Clearly, exactly one of the following sets is in $\mathcal{U}$:
     	\[  \bigcup_{j \in \omega} [n_{2j}, n_{2j+1}) \ \ \text{ and } \ \bigcup_{j \in \omega} [n_{2j+1},n_{2j+2}).  \]
     	Let $U$ denote that element of $\mathcal{U}$. Next, we will define a sequence $(t_i)_{i \in \omega}$ ($t_i \in M_i$) by induction on $i$ so that the following holds
     	\begin{equation} \label{tsor} \begin{array}{rl}
		  i \notin U & \rightarrow \ t_i = r_i, \\
		  i \in U & \rightarrow \ C_{t_0t_1 \dots t_i} \in \mathcal{D}_0^i.
	     \end{array}  \end{equation}
     	Let $t_0 = 1$. Assume that $t_0,t_1, \dots, t_{i-1}$ are already defined. 
     	If $i \notin U$, then let $t_i = r_i$. Otherwise, first recall that by $\eqref{suru}$, for all $C \in \mathfrak{C}^{m_k}$
     	there is an element $C' \in \mathcal{D}^{k+1}_0 \subseteq \mathfrak{C}^{m_{k+1}}$ for which $C' \subseteq C$. Now using $\eqref{frakC}$ and $\eqref{kitolt}$ let
     	$t_i$ be such that $C_{t_0t_1 \dots t_{i-1}t_i} \in \mathcal{D}_0^i$.
     	
     	Now as $(C_{t_0t_1 \dots t_i})_{i \in \omega} \in \prod_{i \in \omega} \mathfrak{C}^{m_i}$, that is a decreasing sequence by $\eqref{kitolt}$, then if $\{ z \} = \bigcap_{i \in \omega} C_{t_0t_1 \dots t_i}$, then $z \in S'$, because by $\eqref{tsor}$
     	\[ \{ i : \ C_{t_0t_1 \dots t_i} \in \mathcal{D}^i_0 \} \supseteq U \in \mathcal {U} .\]
     	But $z \in \varphi_0[R]$ also, since $\omega \setminus U$ contains infinitely many intervals of the form $[n_l, n_{l+1})$, hence
             	$t_{|[n_l,n_{l+1})} = r_{|[n_l,n_{l+1})}$ for infinitely many $l$-s, thus by $\eqref{kockaslemma}$,
     	\[ \{z\} =  \bigcap_{i \in \omega} C_{t_0t_1 \dots t_i} \subseteq \varphi_0[R]. \]
	     This verifies that $\varphi_0[R] \cap S' \neq \emptyset$.

\ep

Now we are ready to state our first main theorem.

\bt \label{atetel}
	If $G$ is a locally compact group that is non-discrete, and $\mu$ is a left Haar measure
	then there is a subgroup $S \leq G$ for which $\mu(S) = 0$, but $S$ is non-meager.
\et
\bp
	Using Lemma $\ref{struct}$, there is an open FL-subgroup $H \leq G$, i.e. for each neighborhood $U\subseteq H$ of the identity, there is a compact normal subgroup $N_U \subseteq U$ of $H$ such that $H/N_U$ is a Lie group having finitely many connected components (fix this operation $U \mapsto N_U$).
	Thus, using Corollary $\ref{kov}$, $H$ is an inverse limit of Lie groups all of which have finitely many connected components.
	Now we define a compact normal subgroup $K \vartriangleleft H$ such that $H/K$ is a second countable Lie group of positive dimension,
	or an inverse limit of a countably infinite sequence of finite groups, where the size of the groups form a strictly increasing sequence, thus (by Proposition $\ref{profinite}$ and $\ref{Lieeset}$) there is a null, but non-meager subgroup in $H/K$.
	Now,  if there is an 
		\[ N \in \mathfrak{N}_H = \{N \vartriangleleft H :\  H/N \text{ is a Lie group } \]
		\[ \text{  with finitely many connected components} \} \] for which $H/N$ is a Lie group of positive dimension, then we are done.
	
	Otherwise, each $H/N$ ($N \in \mathfrak{N}_H$) is a finite group. 
		Recall that, using Lemma $\ref{directed}$ $\mathfrak{N}_H$ is closed under finite intersection.
		First, then there cannot exists a least $N_0 \in \mathfrak{N}_H$ (wrt. the inclusion), because the mapping
		 $g \mapsto (gN)_{N \in \mathfrak{N}_H}$ (which maps $H$ to $\underleftarrow{\lim}_{N \in \mathfrak{N}_H} H/N$) is an isomorpism by Corollary $\ref{kov}$, thus its kernel is $\{ e \} = \bigcap_{N \in \mathfrak{N}_H}N = N_0$ (by Lemma $\ref{reszrsz}$).
		This would imply that $H = H/N_0$ is a finite discrete group. But $H$ was an open subgoup of $G$, thus $G$ is discrete, which is a contradiction.
		
				Then there cannot exists a minimal $N \in \mathfrak{N}_H$ either, because that would be a least element, since $\mathfrak{N}_H$ is closed under finite intersection.
		Therefore, one can construct a strictly decreasing sequence  
		\[ N_0 \gneq N_1 \gneq N_2 \gneq \dots \gneq N_k \gneq \dots \]
		in $\mathfrak{N}_H$, hence each finite group $G/N_{i}$ is a nontrivial homomorphic image of $G/N_{i+1}$ ($i \in \omega$). 
		Then let $K = \bigcap_{i \in \omega} N_i$, and using Lemma $\ref{reszrsz}$ for the group $H$ with $\mathfrak{M} = \{ N_i: \ i \in \omega\}$ we obtain
		\[ H/K \simeq \underleftarrow{\lim}_{i \in \omega} H/N_i .\]
		We have a quotient $H/K$  of $H$ which is an inverse limit of an infinite sequence of finite groups, where the 
		size of the groups is a strictly increasing sequence.
	Let $S' \leq H/K$ be a null but non-meager subgroup (such an $S'$ exists by Proposition $\ref{profinite}$, or $\ref{Lieeset}$).
	Let $S$ be the preimage of $S'$ under the canonical projection $\varphi_K: H \to H/K$ 
	\[ S = \varphi_K^{-1}(S') \leq H .\]
	Using Corollary $\ref{nullpullback}$, $S$ is null in $H$, but then since $H$ is an open subgroup of $G$, $S$ is null in $G$ (the restriction  of a Haar measure of $G$ to an open subgroup $H$ is a Haar measure of $H$, that is unique up to a positive multiplicative constant).
	
	Then we can apply Lemma $\ref{vege}$, i.e. $S$ is not meager in the open subgroup $H$, thus is non-meager in $G$.
\ep

\section{Meager but non-null subgroups}

The following theorem due to Friedman can be found in \cite{Burke} and \cite{Pawl}:
\bt
	In the Cohen model (i.e. adding $\omega_2$ Cohen reals to a model of $ZFC + CH$)
     every $F_\sigma$ subset of $2^\omega \times 2^\omega$ which contains a non-null rectangle must contain a measurable non-null rectangle, i.e.
	\begin{equation} \label{rect} \begin{array}{l} 
	\forall H \subseteq 2^\omega \times 2^\omega, \ \text{ if }\ H \text{ is }F_\sigma: \\
	(\exists C \times D \subseteq H, \ C \times D \notin \mathcal{N}) \Rightarrow (\exists A \times B \subseteq H \text{ measurable}: \ \mu(A \times B) > 0) \\
	
	\end{array} \end{equation}
\et

Next we prove that if this property holds in a locally compact Polish group $G$ then every meager subgroup of $G$ is of measure zero, this lemma is also from \cite{Burke} (stating it only for $G= 2^\omega$, but the proof is the same).

\bl
	Let $G$ be a locally compact  Polish group. Assume that every $F_\sigma$ subset of $G \times G$ which contains a non-null rectangle contains a Haar-measurable non-null rectangle too, i.e.
	\begin{equation} \label{tegla}
			 \begin{array}{l} 
			 \forall H \subseteq G \times G, \ \text{ if } H \text{ is }F_\sigma: \\
			  (\exists C \times D \subseteq H, \ C \times D \notin \mathcal{N}) \Rightarrow (\exists A \times B \subseteq H \text{ measurable}: \ \mu(A \times B) > 0) \\
			 
			 \end{array}
	\end{equation} 
    Then every meager subgroup of $G$ is null.
\el
\br \upshape \label{tegl}
	If $X$ is a Polish space, and $\mu$ is a $\sigma$-finite Borel measure, then for any rectangle $C \times D \subseteq X \times X$,
	it is non-null (wrt. the product measure $\mu \times \mu$) iff $C \notin \mathcal{N}_\mu$ and $D \notin \mathcal{N}_\mu$.
	
\er

\bp
	Let $S \leq G$ be a meager subgroup, assume on the contrary that $S$ is non-null.
	Let $(H_i)_{i \in \omega}$ be nowhere dense closed subsets such that $\bigcup_{i \in \omega} H_i \supseteq S$. Now, if $m: G \times G \to G$ denotes  the multiplication function, then $m^{-1}(\bigcup_{i \in \omega} H_i)$ is an $F_\sigma$ set, containing $S \times S$, which is non-null.
	But then there is a measurable rectangle $A \times B \subseteq m^{-1}(\bigcup_{i \in \omega} H_i)$, which is of positive measure, thus by our Remark $\ref{tegl}$, $\mu(A), \mu(B) > 0$.
	Then, due to a Steinhaus type theorem \cite{bcs}, $AB=m(A \times B)$ has nonempty interior. But $AB \subseteq \bigcup_{i \in \omega} H_i$ which is meager, a contradiction (by the Baire category theorem).
\ep

First (in Lemma $\ref{C-lengyel}$) we show that if every $F_\sigma$ set $H \subseteq 2^\omega \times 2^\omega$ containing a rectangle of positive outer measure contains
a measurable rectangle of positive measure, then this holds for arbitrary locally compact Polish groups. This yields that it is consistent with $ZFC$ that in a locally compact Polish group  meager subgroups are always null.

Later in Theorem $\ref{konzisztens}$ we reduce the general locally compact case to the case of Polish locally compact groups, by showing that if in locally compact Polish groups 
 meager subgroups are null, then this is true in every locally compact group.
\br \upshape \label{merh}
	If $X$ is a Polish space, $\nu$ is a $\sigma$-finite Borel measure on $X$, then for every measurable set $H$ there exists a Borel $B$ such that
	$H \bigtriangleup B$ is null (wrt. $\nu$).
\er

\bl \label{C-lengyel}
Assume that condition ($\ref{tegla}$) holds in $2^\omega$ (i.e. every $F_\sigma$ subset of $2^\omega \times 2^\omega$ which contains a non-null rectangle must contain a measurable Haar-positive rectangle).  Then $\eqref{tegla}$ holds in every non-discrete locally compact Polish group $G$.
\el
\bp
	Let $\mu$ denote the left Haar measure on $G$, and let $\nu$ denote the Haar measure on $2^\omega$.
	Since locally compact Polish groups are $\sigma$-compact, the Haar measure is $\sigma$-finite. 
	
	\bc
	There is a sequence of pairwise disjoint compact sets $C_i$ ($i \in \omega$) in $2^\omega$, for which
	$\nu(2^\omega \setminus \bigcup_{i \in \omega} C_i) = 0$, and similarly a sequence of pairwise disjoint compact sets $K_i$ ($i \in \omega$)
	in $G$ with $\mu(K_i)>0$ and $\mu(G \setminus \bigcup_{i \in \omega} K_i)=0$, and there exist  homeomorphisms $f_i: C_i \to K_i$,
	and positive constants $r_i$ such that
	\[ \forall B \subseteq C_i \text{ Borel }: \ \nu(B) = r_i \mu(f_i(B)). \]
	\ec
	\bp
	First, using the inner regularity, and the $\sigma$-compactness of $G$, let	$E_i$ ($i \in \omega$)  be a sequence of pairwise disjoint compact sets in $G$, for which $\mu(E_j)> 0$, and $\mu(G \setminus \bigcup_{i \in \omega}E_i) = 0$.
	Similarly, let $F_i$ be a sequence of pairwise disjoint compact subsets of $2^\omega$ such that $\nu(2^\omega \setminus \bigcup_{i \in \omega} F_i) = 0$, and $\nu(F_i) > 0$  for all $i$.
	
	Now, since $G$ is non-discrete, every open set is infinite in $G$, and an open set $U$ with compact closure is thus an infinite set.
	But compact sets have finite measure. Then, by the invariance of the measure, every point has the same measure, that must be $0$ because there exists infinite sets  with finite measure, thus we obtain that $\mu$	is a continuous measure. We can apply  the isomorphism theorem for measures \cite[Thm 17.41.]{Kechris}.
	Let $g_i : F_i \to E_i$ be a Borel isomorphism (a bijection which is Borel, and so is its inverse), for which there exists $r_i>0$ such that 
	\[ \nu(H) = r_i \mu(f_i(H)) \ (H \subseteq F_i \text{ Borel}) \]
	Now, by Lusin's theorem (see \cite[Thm 17.12]{Kechris} ) for every measurable set $H \subseteq F_i$ and $\varepsilon>0$, the function $g_i|_H$ is continous on a compact subset
	$H' \subseteq H$  with $\nu(H') > \nu(H)- \varepsilon$.
	Using this, there are  pairwise disjoint compact subsets $(F_i^j)_{j \in \omega}$
	in $F_i$, such that $g_i|_{F_i^j}$ is continuous, and $\nu(F_i \setminus \bigcup_{j \in \omega} F_i^j)=0$.
	Let the sequence $C_i$ ($i \in \omega$) be the enumeration of the $F_i^j$-s in type $\omega$, and let $f_i = g_k^j$, if
	$C_i = F_k^j$, then choosing $K_i$ to be $f_i(C_i)$ works.
	\ep

	 Let $C = \bigcup_{i \in \omega} C_i \subseteq 2^\omega$, $K = \bigcup_{i \in \omega} K_i \subseteq G$,  $f= \bigcup_{i \in \omega} f_i$, then
	 \[ f: C \to K \]
	 is a bijection which is almost everywhere defined in $2^\omega$, and
	 \begin{itemize}
	 	\item for each $S \subseteq C \subseteq  2^\omega$,
	 	if  $\nu(S)=0$, then $\mu(f(S))=0$,
	 	\item if $S$ is measurable, $\nu(S)>0$, then $f(S)$ is measurable,  $\mu(f(S))>0$ (by our Remark $\ref{merh}$, $S$ differs from a Borel by a null set, $f^{-1}$ is a Borel function, and $f$ maps a null set to a null set),
	 	\item for every $F_\sigma$ set $H \subseteq G$, $f^{-1}(H)$ is also $F_\sigma$ (in $2^\omega$).
	 \end{itemize}

	 Now let $H \subseteq G \times G$ be an $F_\sigma$ subset, and $D\times E \subseteq H$  be a rectangle of positive outer measure. Then $H\cap (K\times K)$ is still an $F_\sigma$ set, and by  Remark $\ref{tegl}$, $(D\times E) \cap (K \times K) = (D \cap K) \times 
	 (E \cap K)$ is of positive outer measure, since $\mu(G \setminus K) = 0$. Thus from now on, we can assume that $H \subseteq K \times K$: it suffices to find measurable non-null rectangles in such $H$-s. 
	 But $(f \times f)^{-1}(H)$ is $F_\sigma$,
	 since $f \times f = (\bigcup_{i \in \omega} f_i) \times (\bigcup_{i \in \omega} f_i) = \bigcup_{i,j \in \omega} f_i \times f_j$ where the $f_i$-s are homeomorphisms between compact sets, so are the $f_i \times f_j$-s.
	 And $(f \times f)^{-1}(D \times E) = f^{-1}(D) \times f^{-1}(E)$ is a product of sets of positive outer measure, since $f$ maps a nullset to a nullset. But then $(f \times f)^{-1}(H) \subseteq 2^\omega \times 2^\omega$, that contains  $(f \times f)^{-1}(D \times E)$ contains a measurable 
	 $A \times B$ with $\nu(A) , \nu(B)>0$. Now, using that $f$ maps a positive, measurable set to a positive, measurable set, we obtain that $\mu(f(A)), \mu(f(B)) > 0$, $f(A) \times f(B) \subseteq H$, hence $H$ contains a non-null measurable rectangle, indeed.
\ep

\bt \label{konzisztens}
	If in every non-discrete locally compact Polish group meager subgroups are null, then this is true in non-discrete locally compact groups as well. 
\et
\bp

	Let $S$ be a meager subgroup of $G$. Then, applying Lemma $\ref{struct}$, we get an open FL subgroup $H \leq G$.
	
	 Then, since $H$ is open, $S \cap H \leq H$ is clearly meager (in $H$). 
	 
	 Thus we can apply Lemma $\ref{seged}$ for the co-meager set $R = H \setminus S$ in $H$ and $\mathfrak{J}_0 = \emptyset$,
	 yielding a compact normal subgroup $K' \vartriangleleft H$, and a co-meager set $R' \subseteq H/K'$, for which $H/K'$ is a Polish space, and $\varphi_{K'}^{-1}(R') \subseteq R$, ($\varphi_{K'}: H \to H/K'$). 
	 Thus
	 \[ \varphi^{-1}_{K'}(R') \cap (S \cap H) = \emptyset, \]
	 therefore
	 \[ S \cap H \subseteq H \setminus \varphi^{-1}_{K'}(R') = \varphi^{-1}_{K'}((H/K') \setminus R'),\]
	 thus  projecting $S \cap H$ to $H/K'$ yields a subgroup in $H/K'$ that is disjoint from the co-meager set $R'$.
	 
	We obtain that $\varphi_{K'}(S \cap H) \leq H/K'$ is a meager subgroup in a locally compact Polish group.
	 	 Since in a discrete space a nowhere dense set must be empty, but we found a meager subgroup
	 	  ($e_{H/K'} \in \varphi_{K'}(S \cap H)$), thus $H/K'$ can not be discrete.
	  Then $\varphi_{K'}(S \cap H) \leq H/K'$ is null by our assumptions (wrt. the Haar measure of $H/K'$). 
	 	 But then, by Corollary $\ref{nullpullback}$, the pull-back of $\varphi_{K'}(S \cap H)$ is also null in $H$
	 	 \[ \mu_H( \varphi_{K'}^{-1}(\varphi_{K'}(S \cap H))) = 0,\]
      and $S \cap H \subseteq \varphi_{K'}^{-1}(\varphi_{K'}(S \cap H))$, thus $\mu_H(S \cap H) = 0$.
	
	Since $H$ is an open subgroup of $G$, $S \cap H$ is null with respect to the Haar measure of $G$.
	Let $L \subseteq G$ be a set
	containing exactly one point from each left coset of $H$:
	$G = \bigcup_{l \in L}lH$ (i.e. it is a partition of $G$ into pairwise disjoint open subsets).
	 Now from each left coset $lH$ if $lH \cap S \neq \emptyset$ then pick $s_l \in lH \cap S$. It is easy to see that for such a left coset $lH$ of $H$
	 \[ lH \cap S = s_l\cdot (H \cap S).\]
	  Hence, by the left-invariance of $\mu_G$, $(lH) \cap S$ is null, thus $S$ is locally null.
	  Now, from the fact that the Haar measure is inner regular wrt.  compact sets, we can conclude that $S$ is null in $G$.
	
\ep

\bcor
	It is consistent with $ZFC$ that there are no meager subgroups of positive outer measure in locally compact groups.
\ecor

\br \upshape
	The existence of a meager, non-null subgroup in a locally compact Polish group $G$ implies that there is such a group if extend $G$ by a compact group, formally, if $G \simeq G' / K$, where $K \vartriangleleft G'$, then the preimage of a meager, non-null subgroup under the canonical projection is meager and non-null.  
\er	

\vspace{1cm}

\bex \label{cex}
	In the case of large, non-$\sigma$-compact locally compact groups, requiring outer regularity with respect to open sets, but the inner regularity with respect to compact sets holds only for open sets may result in a different measure than what we used.
	Consider the following example: take the reals with the discrete topology ($\mathbb{R}_d$), and with the euclidean topology ($\mathbb{R}$), then in $\mathbb{R}_d \times \mathbb{R}$:

	$S= R \times \{ 0\}$ is a meager subgroup (it is nowhere dense, since there is a base consisting of sets of the form $\{r\} \times U$, $U\subseteq \mathbb{R}$ open). 
	 But any open $U$ covering of $S$
	can be partitionated into continuum many pairwise disjoint nonempty open sets:
	 $U = \bigcup_{r \in R} U \cap (\{r\} \times \mathbb{R})$, thus $U$ is of infinite measure (recall that open sets have positive measure).
	 Now, from the outer regularity wrt. open sets, $S$ is of infinite measure.
	
\eex

Due to Christensen we can extend the notion of the null ideal to every (not necessarily locally compact) Polish group. The following is still open.
\bq
 What can we say about  small-large subgroups of non-locally compact Polish groups, replacing "null wrt. the Haar measure" by "Haar-null" in the sense of Christensen? Do (always) exist Haar-null but non-meager, and meager but non-Haar-null subgroups in non-locally compact Polish groups?
 (A subset $X$ of a Polish group is Haar-null, if there is a Borel probability measure $\mu$ on $G$, and a Borel set $B \supseteq X$, such that for each $g,h \in G$ $\mu(gBh) = 0$ \cite{Christ}.)
\eq
Since the model of Ros\l{}anowski and Shelah \cite{small-large} was obtained by an  $\omega^\omega$-bounding forcing (thus no Cohen reals added), it is also worth to consider whether the rectangle inclusion property also holds, or at least whether meager subgroups in all locally compact Polish groups are automatically null there.  It would also be interesting to consider whether the non-existence of meager, but non-null subgroups in locally compact Polish groups imply $\eqref{rect}$, the rectangle inclusion property.

\section*{Acknowledgement}

I would like to thank my supervisor Márton Elekes for presenting me the problem, and also giving hints and ideas, providing useful advices.

\end{document}